\documentclass[12pt,a4paper]{article}
\usepackage{bbm}
\usepackage{mathrsfs}
\usepackage{graphicx}
\usepackage{amsmath}
\usepackage{amssymb}
\usepackage{amsfonts}
\usepackage{enumerate}
\usepackage{theorem}

\makeatletter
\def\tank#1{\protected@xdef\@thanks{\@thanks
        \protect\footnotetext[0]{#1}}}
\def\bigfoot{

    \@footnotetext}
\makeatother

\newcommand{\ea}{\end{array}}
\newtheorem{theorem}{Theorem}[section]
\newtheorem{proposition}{Proposition}[section]

\newtheorem{lemma}{Lemma}[section]
\newtheorem{definition}{Definition}[section]
\newtheorem{Rem}{Remark}[section]

{\theorembodyfont{\rmfamily}
}

\title
{\bf The global attractor for the 3-D viscous primitive equations of large-scale moist atmosphere \thanks{This work was partially
supported by NNSF of China(Grant No. 11401057, No. 11301097),   Natural Science Foundation Project of CQ  (Grant No. cstc2016jcyjA0326),
Fundamental Research Funds for the Central Universities(Grant No. 106112015CDJXY100005) , China Scholarship Council (Grant No.201506055003) and GXNSF (Grant No.  2014GXNSFAA118016) .} }
\author{
 Guoli Zhou
\thanks{ Chongqing University, P.R. China }
\tank{E-mail:zhouguoli736@126.com.}
\qquad 
Yanfeng Guo
\thanks{  Guangxi University of  Science and Technology, P.R. China }
\tank{E-mail:guoyan\_feng@aliyun.com.}
}
\begin{document}
\maketitle

\begin{abstract}
Absorbing ball in $H^{1}(\mho)$ is obtained for the strong solution to the three dimensional viscous moist primitive equations under the natural assumption  $Q_{1},Q_{2}\in L^{2}(\mho)$ which is weaker than the assumption $Q_{1},Q_{2}\in H^{1}(\mho)$ in $\cite{GH2}$. In view of the structure of the manifold and the special geometry involved with vertical velocity,  the continuity of the strong solution in $H^{1}(\mho)$ is established with respect to time and initial data. To obtain the existence of the global attractor for the moist primitive equations, the common method is to obtain the absorbing ball in $H^{2}(\mho)$ for the strong solution to the equations. But it is difficult due to the complex structure of the moist primitive equations. To overcome the difficulty, we try to use Aubin-Lions lemma and the continuous property of the strong solutions to the moist primitive equations to prove the the existence of the global attractor which improves the result, the existence of weak attractor, in $\cite{GH2}$.
\end{abstract}


\section{Introduction}
The paper is concerned with the 3-dimensional viscous primitive equations in the pressure coordinate system (see\ e.g. $\cite{GH1, LC, LTW1, LTW2}$ and the references therein).
 \begin{eqnarray}
&&\partial_{t}v+\nabla_{v}v+w\partial_{\xi}v+\frac{f}{R_{0}}v^{\bot} +\mathrm{grad} \Phi+ L_{1}v=0,\\
&&\partial_{\xi}\Phi+\frac{bP}{p}(1+aq)T=0, \\
&&\mathrm{div} v+ \partial_{\xi} w=0, \\
&&\partial_{t} T+\nabla_{v}T +w\partial_{\xi}T-\frac{bP}{p}(1+aq)w-L_{2}T=Q_{1} ,\\
&&\partial_{t}q+ \nabla_{v}q +w\partial_{\xi}q+L_{3}q=Q_{2}  .
\end{eqnarray}
The unknowns for the primitive equations are the fluid velocity field $(v,w )=(v_{\theta},v_{\varphi},w )\in \mathbb{R}^{3}$ with  $v=(v_{\theta},v_{\varphi})$ and $ v^{\perp}=(-v_{\varphi}, v_{\theta} ) $ being horizontal, the temperature $T$, $q$ the mixing ratio of water vapor in the air  and the geopotential $\Phi.$
$f=2 \mathrm{cos\theta}$ is the given Coriolis parameter, $Q_{1}$ corresponds to the sum of the heating of the sun and the heat added or removed by condensation or evaporation, $Q_{2}$ represents the amount of water added or removed by condensation or evaporation, $a$ and $b$ are positive constants with $a\approx 0.618$, $R_{0}$ is the Rossby number, $P$ stands for an approximate value of pressure at the surface of the earth, $p_{0}$ is the pressure of the upper atmosphere with $p_{0}>0$ and the variable $\xi$ satisfies $p=(P-p_{0})\xi+p_{0}$ where $0<p_{0}\leq p\leq P.$ The viscosity, the heat and the water vapor diffusion operators $L_{1},\ L_{2}$ and $L_{3}$ are given respectively as the following:
$$L_{i}=-\nu_{i}\Delta-\mu_{i}\partial_{zz}, i=1,2,3. $$
Here the positive constants $\nu_{1}, \mu_{1}$ are the horizontal and vertical viscosity coefficients;  the positive constant $\nu_{2}, \mu_{2}$ are the horizontal and vertical heat diffusivity coefficients; while the positive constant $\nu_{3}, \mu_{3}$ are the horizontal and vertical water vapor diffusivity coefficients. The definitions of $\nabla_{v}v, \Delta v, \Delta T, \Delta q, \nabla_{v}q, \nabla_{v}T, \mathrm{div} v, \mathrm{grad} \Phi $ will be given in section 2.
\par
The space domain of equations: $(1.1)-(1.5)$ is
$$\mho=S^{2}\times (0,1), $$
 where $S^{2}$ is two-dimensional unit sphere. The boundary value conditions are given by
\begin{eqnarray}
\xi=1 (p=P): \partial_{\xi}v=0,\ \ w=0,\ \ \partial_{\xi}T=\alpha_{s}(T_{s}-T),\ \ \partial_{\xi}q=\beta_{s}(q_{s}-q),
\end{eqnarray}
\begin{eqnarray}
\xi=0 (p=p_{0}): \partial_{\xi}v=0,\ \ w=0,\ \ \partial_{\xi}T=0,\ \ \partial_{\xi}q=0,
\end{eqnarray}
where $\alpha_{s}, \beta_{s}$ are positive constants, $T_{s}$ is the given temperature on the surface of the earth, $q_{s}$ is the given mixing ratio of water vapor on the surface of the earth. To simplify the notations, we set $T_{s}=0$ and $q_{s}=0$ without losing any generality. For the case $T_{s}\neq 0$ and $q_{s}\neq 0,$ we can homogenize the boundary value conditions for $T,q;$ see $\cite{GH1}$ for detailed discussion on this issue. Moreover, using $(1.2), (1.3)$ and the boundary conditions $(1.6)-(1.7)$, we have
\begin{eqnarray}
w(t;\theta,\varphi,\xi)=\int_{\xi}^{1}\mathrm{div}\ v( t;\theta,\varphi,\xi')d\xi',
\end{eqnarray}
\begin{eqnarray}
\int_{0}^{1}\mathrm{div}\ v d\xi=0,
\end{eqnarray}
\begin{eqnarray}
\Phi( t;\theta,\varphi,\xi)=\Phi_{s}( t;\theta,\varphi)+\int_{\xi}^{1}\frac{bP}{p}(1+aq)Td\xi',
\end{eqnarray}
where $\Phi_{s}( t;\theta,\varphi)$ is a certain unknown function at the isobaric surface $\xi=1.$ In this article, we assume that the constants $v_{i}=\mu_{i}=1, i=1,2,3. $ For the general case, the results will still be valid.
Then using $(1.8)-(1.10)$, we obtain the following equivalent formulation for system $(1.1)-(1.7)$ with initial condition
\begin{eqnarray}
\partial_{t}v&+&\nabla_{v}v+ \Big{(} \int_{\xi}^{1}\mathrm{div}\ v( t;\theta,\varphi,\xi')d\xi'\Big{)}  \partial_{\xi}v+\frac{f}{R_{0}}v^{\bot} +\mathrm{grad} \Phi_{s}\nonumber\\
&&+\int_{\xi}^{1}\frac{bP}{p}\mathrm{grad}[(1+aq)T]d\xi'  -\Delta v- \partial_{\xi\xi}v =0,
\end{eqnarray}
\begin{eqnarray}
\partial_{t} T&+&\nabla_{v}T +\Big{(} \int_{\xi}^{1}\mathrm{div}\ v( t;\theta,\varphi,\xi')d\xi'\Big{)}\partial_{\xi}T\nonumber\\
&&-\frac{bP}{p}(1+aq)\Big{(} \int_{\xi}^{1}\mathrm{div}\ v( t;\theta,\varphi,\xi')d\xi'\Big{)}-\Delta T- \partial_{\xi\xi}T =Q_{1} ,
\end{eqnarray}
\begin{eqnarray}
\partial_{t} q+\nabla_{v}q +\Big{(} \int_{\xi}^{1}\mathrm{div}\ v( t;\theta,\varphi,\xi')d\xi'\Big{)}\partial_{\xi}q -\Delta q- \partial_{\xi\xi}q =Q_{2} ,
\end{eqnarray}
\begin{eqnarray}
\int_{0}^{1}\mathrm{div}\ vd\xi=0,
\end{eqnarray}
\begin{eqnarray}
\xi=1 : \partial_{\xi}v=0,\ \ w=0,\ \ \partial_{\xi}T=-\alpha_{s}T,\ \ \partial_{\xi}q=-\beta_{s}q,
\end{eqnarray}
\begin{eqnarray}
\xi=0 : \partial_{\xi}v=0,\ \ w=0,\ \ \partial_{\xi}T=0,\ \ \partial_{\xi}q=0,
\end{eqnarray}
\begin{eqnarray}
v(0;\theta,\phi, \xi)=v_{0}(\theta,\phi, \xi ), T(0;\theta,\phi, \xi)=T_{0}(\theta,\phi, \xi ), q(0;\theta,\phi, \xi)=q_{0}(\theta,\phi, \xi ).
\end{eqnarray}
In order to understand the mechanism of long-term weather prediction, one can take advantage of the historical records and numerical computations to detect the future weather. Alternatively, one should also study the long time behavior mathematically for the equations and models governing the motion. The primitive equations represent the classic model for the study of climate and weather prediction, describing the motion of the atmosphere when the hydrostatic assumption is enforced $\cite{ G, Ha, HW, MT, Ri}.$ But the resulting flow or the atmosphere is rich in its organization and complexity (see $\cite{G, Ha,HW}$), the full governing equations are too complicated to be treatable both from the theoretical and the computational side. To overcome this difficulty, some simple numerical models were introduced. The 2-D and 3-D quasi-geostrophic models have been the subject of analytical mathematical study (see e.g., \cite{BB, C, CMT1, CMT2, CW, EM, M, W1, W2, W3} and references therein).
To the best of our knowledge, the mathematical framework of primitive equations was formulated in $\cite{LTW1, LTW2, LTW3},$
where the definitions of weak and strong solutions were given and the existence of weak solution was proven, leaving the uniqueness of weak solution as an open problem for now. Local well-posedness of strong solutions was obtained in $\cite{GMR, TZ}.$ If the domains was thin, the global well-posedness of 3D primitive equations was shown in $\cite{HTZ}.$  Taking advantage of the fact that the pressure is essentially two-dimensional in the primitive equations, global well-posedness of the full three-dimensional case was established in $\cite{CT1}$ and independently in $\cite{Kob1, Kob2}.$ In the subsequent work $\cite{KZ}$ a different proof was developed which allows one to treat non-rectangular domains. Recently, the results were improved in $\cite{CLT1, CLT2, CLT3, CT2}$  by considering the system with partial dissipation, i.e. , with only partial viscosities or only partial diffusion. For the inviscid primitive equations, finite-time blowup was established in $\cite{CINT}.$ To study the long term behavior of primitive equations, the existence of global attractor was established in $\cite{J}$ and dimensions were proven to be finite in $\cite{JT}.$ When moisture is included, an equation for the conservation of water must be added, which is the case in e.g. $\cite{GH1,GH2, LTW1, PTZ}.$ In $\cite{ZHKTZ},$ global well-posedness of quasi-strong and strong solutions was obtained for the primitive equations of atmosphere in presence of vapour saturation.
\par
The understanding of asymptotic behavior of dynamical system is one of the most important topics of modern mathematical physics. One way to solve the problem for dissipative deterministic dynamical system is to consider its global attractor (see its definition in section 2). Thus, in order to capture the dynamical features of moist primitive equations, Guo and Huang in $\cite{GH2}$ proved the existence of universal attractor which is weakly closed in $V$ (see the notations in section $2$) and attracts any orbit in $V$-weak topology, when time goes to $\infty.$ Then,
one natural question arising from it is the existence of global attractor
which is compact in $V$ and attracts any orbit in $V$-strong topology, as time goes to $\infty.$
\par
As it is stated in $\cite{GH2}$ that laking information for the time derivative of the vertical velocity in the moist primitive equations leads  to the failure of establishing the existence of the global attractor. How to overcome the difficulty?  Inspired by $\cite{J}$,
we try to use Aubin-Lions lemma combined with  continuity of the strong solution in $V$ with respect to time and initial data  to prove the solution operation is compact in $V$. Then the compact property of the solution operator implies that the existence of global attractor, i.e., the universal attractor  is indeed compact in $V$ and attracts all bounded sets in $V$ with respect to $V$-strong topology,  when $t\rightarrow \infty$.

\par
In this  article, we first try to obtain time-uniform $a\  priori$ estimates in various function spaces under the natural assumption $Q_{1},Q_{2}\in L^{2}(\mho),$ reducing the stronger assumption $Q_{1},Q_{2}\in H^{1}(\mho)$ in $\cite{GH2}$. Then, by making delicate and careful estimates of $L^{4}$ norm for the strong solution, estimates about $L^{3}$ norm which is required in $\cite{GH2}$ is omitted. To obtain uniform boundedness with respect to time $t$ in $V,$ estimates about $\partial_{\xi}T, \partial_{\xi}q$ are carefully considered, which is more complex than oceanic primitive equations. We recall that the continuity of the strong solution with respect to initial data was shown to be true in weak solution space $H$ in previous work due to no boundedness is available for the derivatives of the vertical velocity. This is not sufficient for our purpose. To overcome the difficulty, the structure of the manifold and the special geometry involved with the vertical velocity are used to obtain the continuity of the strong solution with respect to initial data in $V$, improving the results before. Finally, in order to prove the absorbing ball is compact in $V,$ the common method is to show that the ball is uniformly bounded with respect to time $t$ in $H^{2}(\mho).$ But it is difficult to achieve because of the high nonlinearity of the moist primitive equations. Inspired by $\cite{J}$, we try to use an Aubin-Lions compactness lemma combined with a continuity argument to show that the solution operator is compact in $V$ for every time $t>0$, which further implies  the existence of the global attractor for the dynamical system generated by the primitive equations of large-scale moist atmosphere.

\par
The remaining of the paper is organized as follows. In section $2$, we present the notations and recall some important facts which are crucial to later analysis. Absorbing ball is obtained in section $3$, where we correct the mistake $fv^{\bot}\times \Delta v=0$ in $\cite{GH2}.$ Section $4$ and section $5$ are for continuity of the strong solution with respect to time $t$ and initial data in $V$ respectively.
Finally, in section $6$, using Aubin-Lions lemma and the continuity properties of the strong solution we prove the existence of global attractor.  As usual, the positive constants $c$ may change from one line to the next, unless, we give a special declaration.

\section{Preliminaries }
In this section we collect some preliminary results that will be used in the rest of this paper, and we start with the following notations which will be used throughout this work.
Denote
\begin{eqnarray*}
\bar{v}=\int_{0}^{1}vd\xi ,\ \   \tilde{v}=v-\bar{v} .
\end{eqnarray*}
Then we have
\begin{eqnarray}
\nabla \cdot  \bar{v}=0,\ \   \bar{\tilde{v}}=0.
\end{eqnarray}
\par
Now we give the definitions of some differential operators. Firstly,  the natural generalization of the directional derivative on the Euclidean space to the covariant derivative on $S^{2}$ is given as follows.  Let $T,q, \in C^{\infty}(\mho),   \Phi_{s} \in C^{\infty}(S^{2} )$ and
\begin{eqnarray*}
v=v_{\theta}e_{\theta}+v_{\varphi}e_{\varphi},\ \ \ \ \ u=u_{\theta}e_{\theta} +u_{\varphi}e_{\varphi}\ \ \ \ \in C^{\infty}(T\mho| TS^{2}),
\end{eqnarray*}
where $C^{\infty}(T\mho| TS^{2}) $ is the first two components of smooth vector fields on $\mho.$  We define  the covariant derivative of $u,T $ and $q$ with respect to $v$  as follows
\begin{eqnarray*}
\nabla_{v}u=(v_{\theta}\partial_{\theta}u_{\theta}+\frac{v_{\varphi}}{\mathrm{sin}\theta}\partial_{\varphi}u_{\theta}-v_{\varphi}u_{\varphi}\mathrm{cot}\theta )e_{\theta}+(v_{\theta}\partial_{\theta}u_{\varphi}+\frac{v_{\varphi}}{\mathrm{sin}\theta}\partial_{\varphi}u_{\varphi}+v_{\varphi}u_{\theta}\mathrm{cot}\theta )e_{\varphi},
\end{eqnarray*}
\begin{eqnarray*}
\nabla_{v}T= v_{\theta}\partial_{\theta}T+ \frac{v_{\varphi}}{\mathrm{sin}\theta}\partial_{\varphi}T,
\end{eqnarray*}
\begin{eqnarray*}
\nabla_{v}q= v_{\theta}\partial_{\theta}q+ \frac{v_{\varphi}}{\mathrm{sin}\theta}\partial_{\varphi}q.
\end{eqnarray*}
We give the definition of the horizontal gradient $\nabla = \mathrm{grad}$ for $T$ and $\Phi_{s}$ on $S^{2}$ by
\begin{eqnarray*}
\nabla T= \mathrm{grad} T= (\partial_{\theta}T) e_{\theta} +\frac{1}{\mathrm{sin}\theta}(\partial_{\varphi}T)e_{\varphi},
\end{eqnarray*}
\begin{eqnarray*}
\nabla \Phi_{s}= \mathrm{grad} \Phi_{s}=( \partial_{\theta}\Phi_{s} )e_{\theta} +\frac{1}{\mathrm{sin}\theta}(\partial_{\varphi}\Phi_{s})e_{\varphi}.
\end{eqnarray*}
We define the divergence of $v$ by
\begin{eqnarray*}
\mathrm{div} v= \mathrm{div} (v_{\theta}e_{\theta}+v_{\varphi}e_{\varphi} )=\frac{1}{\mathrm{sin}\theta}(\partial_{\theta}(v_{\theta}\mathrm{sin} \theta)+\partial_{\varphi} v_{\varphi}  ).
\end{eqnarray*}
The horizontal Laplace-Beltrami operator of scalar functions $T$ and $q$ are
\begin{eqnarray*}
\Delta T= \mathrm{div}(\mathrm{grad}T)=\frac{1}{\mathrm{sin}\theta}[\partial_{\theta}(\mathrm{sin}\theta\partial_{\theta}T )+\frac{1}{\mathrm{sin}\theta}\partial_{\varphi\varphi}T],
\end{eqnarray*}
\begin{eqnarray*}
\Delta q= \mathrm{div}(\mathrm{grad}q)=\frac{1}{\mathrm{sin}\theta}[\partial_{\theta}(\mathrm{sin}\theta\partial_{\theta}q )+\frac{1}{\mathrm{sin}\theta}\partial_{\varphi\varphi}q].
\end{eqnarray*}
We define the horizontal Laplace-Beltrami operator $\Delta$ for vector functions on $S^{2}$ as
\begin{eqnarray*}
\Delta v=(\Delta v_{\theta}-\frac{2\mathrm{cos}\theta}{\mathrm{sin}^{2}\theta}\partial_{\varphi}v_{\varphi}-\frac{v_{\theta}}{\mathrm{sin}^{2}\theta  }  )e_{\theta}+(\Delta v_{\varphi}+ \frac{2\mathrm{cos}\theta}{\mathrm{sin}^{2}\theta}\partial_{\varphi}v_{\theta}-\frac{v_{\varphi}}{\mathrm{sin}^{2}\theta   } )e_{\varphi}.
\end{eqnarray*}
Consequently, by integration by parts, we have
\begin{eqnarray}
\int_{0}^{1}w\partial_{\xi}vd\xi= \int_{0}^{1}v \mathrm{div} v d\xi= \int_{0}^{1}\tilde{v} \mathrm{div} \tilde{v}d\xi,
\end{eqnarray}
\begin{eqnarray}
\int_{0}^{1}\nabla_{v}vd\xi= \int_{0}^{1}\nabla_{\tilde{v}}\tilde{v}d\xi+\nabla_{\bar{v}}\bar{v}.
\end{eqnarray}
Taking the average of equations $(1.11)$ in the $z$ direction, over the interval $(0,1)$ and using $(2.18)-(2.20)$ and the boundary conditions $(1.15 )-(1.16),$  we arrive at
\begin{eqnarray}
\partial_{t}\bar{v}&+&\nabla_{\bar{v}}\bar{v}+\overline{\tilde{v}\mathrm{div}\tilde{v} +  \nabla_{\tilde{v}}\tilde{v}}+\frac{f}{R_{0}} \bar{v}^{\bot}+\mathrm{grad } \Phi_{s}+\int_{0}^{1}\int_{\xi}^{1}\frac{bP}{p}\mathrm{grad}[(1+aq)T ]d\xi'd\xi\nonumber\\
&&-\Delta \bar{v}=0\ \ \mathrm{in}\ S^{2}.
\end{eqnarray}
By subtracting $(2.21)$ from $(1.11),$ we obtain the following equation
\begin{eqnarray}
\partial_{t}\tilde{v}&+&\nabla_{\tilde{v}} \tilde{v}+\Big{(}\int_{\xi}^{1} \mathrm{div} \tilde{v} d\xi' \Big{)}\partial_{\xi} \tilde{v}
+\nabla_{\tilde{v}}\bar{v}+\nabla_{\bar{v}}\tilde{v}-\overline{(\tilde{v}  \mathrm{div} \tilde{v} + \nabla_{\tilde{v}}\tilde{v}  )}
+\frac{f}{R_{0}}\tilde{v}^{\bot}\nonumber\\
&&+\int_{\xi}^{1}\frac{bP}{p} \mathrm{grad}[(1+aq)T  ]d\xi'- \int_{0}^{1}\int_{\xi}^{1}\frac{bP}{p} \mathrm{grad}[(1+aq)T  ]d\xi'd\xi\nonumber\\
&&-\Delta \tilde{v}-\partial_{\xi\xi}\tilde{v}=0\ \ \mathrm{in}\ \mho,
\end{eqnarray}
with the following boundary value conditions
\begin{eqnarray}
\partial_{\xi} \tilde{v}=0\ \mathrm{on}\ \xi=1\ \mathrm{and}\  \xi=0.
\end{eqnarray}
Let $e_{\theta}, e_{\varphi}, e_{\xi}$ be the unite vectors in $\theta, \varphi$ and $\xi$ directions of the space domain $\mho$ respectively,
\begin{eqnarray*}
e_{\theta}=\partial_{\theta},\ \ \ \ e_{\varphi}= \frac{1}{\mathrm{sin\theta}}\partial_{\varphi},\ \ \ \ e_{\xi}=\partial_{\xi}.
\end{eqnarray*}
The inner product and norm on $T_{(\theta, \varphi, \xi)}\mho$  (the tangent space of $\mho$ at the point $(\theta,  \varphi, \xi )$  ) are defined by
\begin{eqnarray*}
(u,v)=u\cdot v=\sum_{i=1}^{3} u_{i}v_{i},\ \ \ \ |u|=(u,\ u)^{\frac{1}{2}},
\end{eqnarray*}
where $u=u_{1}e_{\theta}+u_{2}e_{\varphi}+u_{3}e_{\xi}\in T_{(\theta, \varphi, \xi)}\mho$ and  $v=v_{1}e_{\theta}+v_{2}e_{\varphi}+v_{3}e_{\xi}\in T_{(\theta, \varphi, \xi)}\mho.$ For $1\leq p\leq \infty,$ let $L^{p}(\mho), L^{p}(S^{2})$ be the usual Lebesgue spaces with the norm $|\cdot|_{p}$ and $|\cdot|_{L^{p}(S^{2})}$
respectively. If there is no confusion, we will write $|\cdot|_{p}$ instead of $|\cdot|_{L^{p}(S^{2})}$.  $L^{2}(T\Omega|TS^{2} ) $ is the first two components of $L^{2}$ vector fields on $\mho$ with the norm $|v|_{2}=(\int_{\mho}(|v_{\theta}|^{2}+|v_{\varphi}|^{2} ) d\mho )^{\frac{1}{2}},$  where $v=(v_{\theta}, v_{\varphi} ): \mho\rightarrow TS^{2}$. Denoted by $C^{\infty}(S^{2})$ the functions of all smooth functions from $S^{2}$ to $\mathbb{R}.$ Similarly, we can define $C^{\infty}(\mho)$. $H^{m}(\mho)$ is the Sobolev space of functions which are in $L^{2},$ together with all their covariant derivatives with respect to $e_{\theta}, e_{\varphi}, e_{\xi}$ of order $\leq m,$ with the norm
\begin{eqnarray*}
\|h\|_{m}=[\int_{\mho}(\sum_{1\leq k\leq m}\sum_{i_{j}=1,2,3;j=1,...,k}|\nabla_{i_{1}}\cdots \nabla_{i_{k}}h |^{2}+|h|^{2}   )  ]^{\frac{1}{2}},
\end{eqnarray*}
where $\nabla_{1}=\nabla_{e_{\theta}}, \nabla_{2}=\nabla_{e_{\varphi}}$ and $\nabla_{3}=\partial_{\xi} $ which are defined above. Denote $H^{m}(T\mho|TS^{2}) =\{v;v=(v_{\theta}, v_{\varphi}):\mho\rightarrow TS^{2}, \|v\|_{m}< \infty   \},$ where the norm  is similar to that of $H^{m}(\mho)$ (i.e., let $h=(v_{\theta}, v_{\varphi} )= v_{\theta}e_{\theta}+v_{\varphi}e_{\varphi}.$)

We will conduct our work in the following functional spaces. Let
\begin{eqnarray*}
\mathcal{V}_{1}:=\{v; v\in C^{\infty}(T\mho|TS^{2}),\ \partial_{\xi}v|_{\xi=0}=0,\ \partial_{\xi}v|_{\xi=1}=0,\ \int_{0}^{1}\mathrm{div}\ vd\xi=0  \},
\end{eqnarray*}
\begin{eqnarray*}
\mathcal{V}_{2}:=\{T; T\in C^{\infty}(\mho),\ \partial_{\xi}T|_{\xi=0}=0,\ \partial_{\xi}T|_{\xi=1}=-\alpha_{s}T  \},
\end{eqnarray*}
\begin{eqnarray*}
\mathcal{V}_{3}:=\{q; q\in C^{\infty}(\mho),\ \partial_{\xi}q|_{\xi=0}=0,\ \partial_{\xi}q|_{\xi=1}=-\beta_{s}q  \}.
\end{eqnarray*}
We denote by $V_{1}, V_{2}$ and $V_{3}$ the closure spaces of $\mathcal{V}_{1},  \mathcal{V}_{2}$ and $\mathcal{V}_{3}$ in $H^{1}(\mho)$ under $H^{1}-$ topology, respectively. In addition, we denote by $H_{1}, H_{2}$ and $ H_{3}$ the closure of $\mathcal{V}_{1}, \mathcal{V}_{2}   $  and $\mathcal{V}_{3} $  in $L^{2}(\mho)$  under $L^{2}-$ topology. Let $H:= H_{1}\times H_{2}\times H_{3}$  and $V=V_{1}\times V_{2}\times V_{3}$ with $V'$ being  dual space of $V$.
 By definition, the inner products and norms on $V_{1}, V_{2}$ and $ V_{3}$ are given by
\begin{eqnarray*}
\langle v,v_{1}      \rangle_{V_1}=\int_{\mho}(\nabla_{e_{\theta}}v\cdot \nabla_{e_{\theta}}v_{1}+\nabla_{e_{\varphi}}v\cdot \nabla_{e_{\varphi}}v_{1}+\partial_{\xi}v\partial_{\xi}v_{1}+v\cdot v_{1} )d\mho,
\end{eqnarray*}
\begin{eqnarray*}
\|v\|_{1}=\langle v,     v\rangle_{V_{1}}^{\frac{1}{2}},\ \ \ \ \forall\ v,\ v_{1}\in V_{1},
\end{eqnarray*}
\begin{eqnarray*}
\langle T,   T_{1} \rangle_{V_{2}}=\int_{\mho}(\mathrm{grad} T\cdot \mathrm{grad} T_{1} +\partial_{\xi}T \partial_{\xi}T_{1} )d\mho +\alpha_{s}\int_{S^{2}}TT_{1} dS^{2},
\end{eqnarray*}
\begin{eqnarray*}
\|T\|_{1}=\langle T,   T \rangle_{V_{2}}^{\frac{1}{2}},\ \ \ \forall\ T,T_{1}\in V_{2},
\end{eqnarray*}
\begin{eqnarray*}
\langle q,   q_{1} \rangle_{V_{3}}=\int_{\mho}(\mathrm{grad} q\cdot \mathrm{grad} q_{1} +\partial_{\xi}q \partial_{\xi}q_{1})d\mho + \beta_{s}\int_{S^{2}}qq_{1}dS^{2} ,
\end{eqnarray*}
\begin{eqnarray*}
\|q\|_{1}=\langle q,   q\rangle_{V_{3}}^{\frac{1}{2}},\ \ \ \forall\ q,q_{1}\in V_{3}.
\end{eqnarray*}
  Let $V_{i}'(i=1,2,3)$ be the dual space of $V_{i}$ with $\langle , \rangle$ being the inner products between $V_{i}'$ and $V_{i}.$ Without confusion, we also  denote by $\langle, \rangle$ the inner product in $L^{2}(\mho)$ and $L^{2}(S^{2}).$  Define the linear operator $A_{i}: V_{i}\mapsto V_{i}', i=1,2,3$ :
\begin{eqnarray*}
&&\langle A_{1}u_{1}, u_{2}    \rangle= \langle u,  v  \rangle_{V_{1}},\ \ \ \forall\ u_{1},\ u_{2}\in V_{1};\\
&&\langle A_{2}\theta_{1} ,  \theta_{2}   \rangle= \langle \theta_{1},  \theta_{2}   \rangle_{V_{2}},\ \ \ \forall\ \theta_{1},\theta_{2} \in V_{2};\\
&&\langle A_{3} q_{1} ,  q_{2}   \rangle= \langle q_{1},  q_{2}   \rangle_{V_{3}},\ \ \ \forall\ q_{1},q_{2} \in V_{3}.
\end{eqnarray*}
Denote $D(A_{i})=\{\eta\in V_{i},    A_{i}\eta\in H_{i} \}.$ Since $ A_{i}$ is positive self-adjoint with compact resolvent, according to the classic spectral theory we can define the power $A_{i}^{s}$ for any $s\in \mathbb{R}.$ Then we have $D(A_{i}^{\frac{1}{2}})=V_{i}$ and $D(A_{i}^{-\frac{1}{2}})=V_{i}'.$ Moreover,
\begin{eqnarray*}
D(A_{i})\subset V_{i}\subset H_{i}\subset V_{i}'\subset D(A_{i})',
\end{eqnarray*}
where $D(A_{i})' $ is the dual space of $ D(A_{i})$ and the embeddings above are all compact.
In the following, we state some lemmas including integrations by parts and uniform Gronwall lemma, which are frequently used in our paper. For the proof of Lemma $2.1$-Lemma $2.3$, we can see $\cite {GH1}.$ The proof of uniform Gronwall lemma was given in $\cite{FP, T}.$
\begin{lemma}
Let $u=(u_{\theta}, u_{\varphi} ), v= (v_{\theta}, v_{\varphi})\in C^{\infty}(T\mho|TS^{2} )$ and $p\in C^{\infty}(S^{2}).$ Then
\begin{eqnarray*}
\int_{S^{2}}p\ \mathrm{div}\ udS^{2}=-\int_{S^{2}}\nabla p \cdot u dS^{2},
\end{eqnarray*}
\begin{eqnarray*}
\int_{\mho}\nabla p \cdot v d\mho=0\ \ \ for\ any\ v\in V_{1},
\end{eqnarray*}
and
\begin{eqnarray*}
\int_{\mho}(-\Delta u) \cdot vd\mho= \int_{\mho}(\nabla_{e_{\theta} } u\cdot \nabla_{e_{\theta}}v+ \nabla_{e_{\varphi} } u\cdot \nabla_{e_{\varphi}}v+u\cdot v    )d\mho.
\end{eqnarray*}
\end{lemma}
\begin{lemma}
For any $h\in C^{\infty}(S^{2}), v\in C^{\infty}(T\mho|TS^{2}),$ we have
\begin{eqnarray*}
\int_{S^{2}}\nabla_{v}hd S^{2}+\int_{S^{2}}h\ \mathrm{div}\ v dS^{2} = \int_{S^{2}}\mathrm{div} (hv)dS^{2}=0.
\end{eqnarray*}
\end{lemma}
\begin{lemma}
Let $u,v \in V_{1}, T\in V_{2}, q\in V_{3}.$ Then we have
\begin{eqnarray*}
\int_{\mho}[\nabla_{u}v+(\int_{\xi}^{1} \mathrm{div }\ u d\xi'  )\partial_{\xi}v ]v d\mho=0,
\end{eqnarray*}
\begin{eqnarray*}
\int_{\mho}[\nabla_{u}g+(\int_{\xi}^{1} \mathrm{div}\ u d\xi'  )\partial_{\xi}g ]g d\mho=0,\ \ for\ g=T\ or\ g=q,
\end{eqnarray*}
\begin{eqnarray*}
\int_{\mho}\Big{(} \int_{\xi}^{1}\frac{bP}{p} \mathrm{grad} [(1+aq)T]d\xi'\cdot u-\frac{bP}{p}(1+aq)T(\int_{\xi}^{1} \mathrm{div }\ u d\xi'  )   \Big{)}=0.
\end{eqnarray*}
\end{lemma}
\begin{lemma}
 Let $f,g$ and $h$ be three non-negative locally integrable functions on $(t_{0}, \infty)$ such that
\begin{eqnarray*}
 \frac{df}{dt}\leq gf +h,\ \ \ \forall\ t\geq t_{0},
\end{eqnarray*}
and
\begin{eqnarray*}
\int_{t}^{t+r}f(s)ds\leq a_{1},\ \ \int_{t}^{t+r}g(s)ds\leq a_{2},\ \ \int_{t}^{t+r}h(s)ds\leq a_{3},\ \ \forall\ t\geq t_{0},
\end{eqnarray*}
where $r, a_{1}, a_{2}, a_{3}$ are positive constants. Then
\begin{eqnarray*}
 f(t+r)\leq (\frac{a_{1}}{r}+a_{3} )e^{a_{2}},\ \ \ \forall\ t\geq t_{0}.
\end{eqnarray*}
\end{lemma}
Before considering the long time behavior of the dynamics, we recall the definitions of strong solution to $(1.11)-(1.17).$
\begin{definition}
Suppose $Q_{1},Q_{2}\in L^{2}(\mho), (v_{0}, T_{0}, q_{0})\in V$ and $\tau>0.$ $(v,T, q)$ is called a strong solution of $(1.11)-(1.17)$ on the time interval $[0, \tau]$ if it satisfies $(1.11)-(1.13)$ in a week sense, and also
\begin{eqnarray*}
v\in C([0, \tau]; V_{1})\cap L^{2}([0,\tau]; H^{2}(\mho)),
\end{eqnarray*}
\begin{eqnarray*}
T\in C([0, \tau]; V_{2})\cap L^{2}([0,\tau]; H^{2}(\mho)),
\end{eqnarray*}
\begin{eqnarray*}
q\in C([0, \tau]; V_{3})\cap L^{2}([0,\tau]; H^{2}(\mho)),
\end{eqnarray*}
\begin{eqnarray*}
\partial_{t} v, \partial_{t} T,  \partial_{t} q \in L^{1}([0, \tau]; L^{2}(\mho)).
\end{eqnarray*}
\end{definition}
Now we state the global well-posedness theorem for the strong solution as follows. For the proof of the therem, one can refer to $\cite{GH2}.$
\begin{proposition}
Let $Q_{1}, Q_{2} \in H^{1}(\mho), U_{0}=(v_{0}, T_{0}, q_{0}) \in V .$
Then for any $T > 0$ given, the global strong solution $U$ of the system $(1.11)-(1.17)$ is unique on the interval $[0, T]$. Moreover, the strong solution $U$ is continuous with respect to initial data in $H.$
\end{proposition}
\begin{Rem}
Notice that there are some gap between the Proposition $2.1$ and Definition 2.1. In fact, in section 3 of this work we will find that if the condition  $Q_{1}, Q_{2}\in H^{1}(\mho)$  is relaxed as $Q_{1}, Q_{2}\in L^{2}(\mho) $, the result of Proposition $2.1$ still holds.
\end{Rem}
\begin{Rem}
In the Proposition $2.1$, the result that the strong solution is continuous with respect to initial data in $H$ is not sufficient for our purpose. We will improve the result in section 5 by establishing the strong solution is continuous with respect to initial data in $V.$
\end{Rem}
To see the difference between global attractor and universal attractor, we introduce the two definitions in the following. For more details, we refer to $\cite{GH2, H, T}$ and other references.
Let $(X,d)$ be a separable metric space and $ S(t): X\rightarrow X, 0\leq t<\infty,$ be a semigroup satisfying:
\par
$(i)\ S(t)S(s)x=S(t+s)x, $ for all $t,s \in \mathbb{R}_+$ and $x\in X;$
\par
$(ii)\ S(0)=I$ (Identity in $X$);
\par
$(iii)\ S(t)$ is continuous in $X$ for all $t\geq0 .$
\par
Typically, $S(t)$ is associated with a autonomous differential equation; $S(t)x$ is the state at time $t$ of the solution whose initial data is $x.$

\begin{definition}
   A subset $\mathcal{A}$ in $X$ is said to be a global attractor if it satisfies the following properties:\\
$(i)$ $\mathcal{A}$ is compact in $X;$\\
$(ii)$ for every $t\geq 0, S(t) \mathcal{A} =\mathcal{A};$\\
$(iii)$ for every bounded set $B$ in $X,$ the set $S(t)B$ converges to $\mathcal{A}$ in $X$, when $t\rightarrow \infty, i.e.,$
$$\lim\limits_{t\rightarrow \infty}d(S(t)B, \mathcal{A}) =0. $$
Here, and in the following, for $A$and $ B$ subsets of $X, d(A, B)$ is the semi-distance given by
$$d(A,B)=\sup\limits_{x\in A}\inf\limits_{y\in B} d(x,y).$$
\end{definition}

\begin{definition}
   A subset $\mathcal{A}$ in $X$ is said to be a universal attractor or weak attractor if it satisfies the following properties:\\
$(i)$ $\mathcal{A}$ is bounded and weakly closed in $X;$\\
$(ii)$ for every $t\geq 0, S(t) \mathcal{A} =\mathcal{A};$\\
$(iii)$ for every bounded set $B$ in $X,$ the set $S(t)B$ converges to $\mathcal{A}$ with respect to $X$-weak topology, when $t\rightarrow \infty, i.e.,$
$$\lim\limits_{t\rightarrow \infty}d_{X}^{w}(S(t)B, \mathcal{A}) =0 $$
where the distance $d_{X}^{w} $ is induced by the $X$-weak topology.
\end{definition}

\section{Uniform estimates and absorbing balls}
In this section, we will obtain some useful uniform $a$ $priori$ estimates about the solution to $(1.11)-(1.17)$ under the natural assumption $Q_{1}, Q_{2}\in L^{2}(\mho)$ and give a proof of the existence of the absorbing ball in $V$ for the solution to the moist primitive equation. The estimates of this section are rigorous without justification using Galerkin approximation due to Proposition $2.1.$
\subsection{$L^{2}$ estimates of $v, T, q$}
Taking inner product of $(1.13)$ with $q$ in $L^{2}(\mho),$  by Lemma $2.3$ we have
\begin{eqnarray}
\frac{1}{2}\frac{d |q|_{2}^{2}}{dt}+|\nabla q|_{2}^{2}+|\partial_{\xi}q|_{2}^{2}+\beta_{s}|q|_{\xi=1}|_{2}^{2}=\int_{\mho}qQ_{2}d\mho.
\end{eqnarray}
Since $q(\theta,\varphi,\xi)=-\int^{1}_{\xi}\partial_{\xi}qd\xi'+q|_{\xi=1},$ by H$\mathrm{\ddot{o}}$lder inequality and Cauchy-Schwarz inequality we have
\begin{eqnarray}
|q|_{2}^{2}\leq 2|\partial_{\xi} q|_{2}^{2}+2|q_{\xi=1}|_{2}^{2},
\end{eqnarray}
which together with $(3.24)$ implies that there exists a positive constant $c$ such that
\begin{eqnarray*}
\frac{d}{dt}|q|_{2}^{2}+c|q|_{2}^{2}\leq |Q_{2}|_{2}^{2}.
\end{eqnarray*}
Therefore, we have
\begin{eqnarray}
|q(t)|_{2}^{2}\leq e^{-ct}|q_{0}|_{2}^{2}+c|Q_{2}|_{2}^{2}.
\end{eqnarray}
By $(3.24)$ and $(3.26),$ for arbitrary $t_{0}\geq 0$ we have
\begin{eqnarray}
&&\int_{t_{0}}^{t_{0}+1}(|\nabla q(t)|_{2}^{2}+|\partial_{\xi}q(t)|_{2}^{2}+\beta_{s}|q|_{\xi=1}(t)|_{2}^{2} )dt\nonumber\\
&&\leq |q(t_{0})|_{2}^{2}+|Q|_{2}^{2}\leq e^{-ct_{0}}|q_{0}|_{2}^{2}+c|Q_{2}|_{2}^{2}.
\end{eqnarray}
Taking an analogous argument as $(3.24)$  , we have
\begin{eqnarray}
\frac{1}{2}\frac{d|T|^{2}_{2}}{dt}+|\nabla T|_{2}^{2}+|\partial_{\xi} T|_{2}^{2}+\alpha_{s}|T|_{\xi=1}|_{2}^{2}=\int_{\mho}\frac{bP}{p}(1+aq)Twd\mho+\int_{\mho}Q_{1}Td\mho.
\end{eqnarray}
Multiplying $v$ with respect to $(1.11)$ and integrating on $\mho,$ by Lemma $2.1,$   Lemma $2.3$ and $(\frac{f}{R_{0}}\times v )\cdot v =0$ we have
\begin{eqnarray}
\frac{1}{2}\frac{d|v|_{2}^{2}}{dt}+|\nabla_{e_{\theta}}v|_{2}^{2}+|\nabla_{e_{\varphi}}v|_{2}^{2}+|v|_{2}^{2}+|\partial_{\xi}v|_{2}^{2}
=-\int_{\mho}\Big{(} \int^{1}_{\xi}\frac{bP}{p} \mathrm{grad}[(1+aq)T ]d\xi'    \Big{)}\cdot v.
\end{eqnarray}
Combining $(3.28)-(3.29)$ and Lemma $2.3$ yields
\begin{eqnarray}
\frac{1}{2}\frac{d(|v|_{2}^{2}+|T|_{2}^{2})}{dt}&+&|\nabla_{e_{\theta}}v|_{2}^{2}+\nabla_{e_{\varphi}}v|_{2}^{2}+|\partial_{\xi}v|_{2}^{2}+|v|_{2}^{2}\nonumber\\
&&+|\nabla T|_{2}^{2}+|\partial_{\xi}T|_{2}^{2}+\alpha_{s}|T|_{\xi=1}|_{2}^{2}=\int_{\mho}Q_{1}Td \mho.
\end{eqnarray}
Similarly to the deduction of $(3.25),$ we have
\begin{eqnarray*}
|T|_{2}^{2}\leq 2|\partial_{\xi} T|_{2}^{2}+2|T_{\xi=1}|_{2}^{2},
\end{eqnarray*}
which together with $(3.29)$ and H$\mathrm{\ddot{o}}$lder inequality implies
\begin{eqnarray}
\frac{d(|v|_{2}^{2}+|T|_{2}^{2})}{dt}&+&|\nabla_{e_{\theta}}v|_{2}^{2}+|\nabla_{e_{\varphi}}v|_{2}^{2}+|\partial_{\xi}v|_{2}^{2}+|v|_{2}^{2}\nonumber\\
&&+|\nabla T|_{2}^{2}+|\partial_{\xi}T|_{2}^{2}+\alpha_{s}|T|_{\xi=1}|_{2}^{2}\leq c|Q_{1}|_{2}^{2}.
\end{eqnarray}
Therefore, we conclude that
\begin{eqnarray}
|v(t)|_{2}^{2}+|T(t)|_{2}^{2}\leq e^{-ct}(|v_{0}|_{2}^{2}+|T_{0}|_{2}^{2})+c|Q_{1}|_{2}^{2}.
\end{eqnarray}
Combining $(3.30)$ and $(3.31)$ we arrive at
\begin{eqnarray}
&&\int_{t_{0}}^{t_{0}+1}|\nabla_{e_{\theta}}v(t)|_{2}^{2}+\nabla_{e_{\varphi}}v(t)|_{2}^{2}+|\partial_{\xi}v(t)|_{2}^{2}+|v(t)|_{2}^{2}dt\nonumber\\
&&+\int_{t_{0}}^{t_{0}+1}|\nabla T(t)|_{2}^{2}+|\partial_{\xi}T(t)|_{2}^{2}+\alpha_{s}|T|_{\xi=1}(t)|_{2}^{2}dt\leq e^{-ct}(|v_{0}|_{2}^{2}+|T_{0}|_{2}^{2})+c|Q_{1}|_{2}^{2}.
\end{eqnarray}
\subsection{$L^{4}$ estimates of $q$}
Multiplying $q^{3}$ on both sides of $(1.13)$ and integrating on $\mho$ yields
\begin{eqnarray}
&&\frac{1}{4}\frac{d |q|_{4}^{4}}{dt}+3||\nabla q|q  |_{2}^{2}+3||\partial_{\xi}q|q|_{2}^{2}+\beta_{s}| q|_{\xi=1}|_{4}^{4}\nonumber\\
&=&\int_{\mho}Q_{2}q^{3}d\mho-\int_{\mho}[\nabla_{v}q+(\int_{\xi}^{1}\mathrm{div} v d\xi' ) \partial_{\xi} q]q^{3}.
\end{eqnarray}
Since by Lemma 2.2 or Lemma 2.3 we have
\begin{eqnarray}
\int_{\mho}[\nabla_{v}q+(\int_{\xi}^{1}\mathrm{div} v d\xi' ) \partial_{\xi}q ]q^{3}=0.
\end{eqnarray}
Then $(3.34)$ and $(3.35)$ imply
\begin{eqnarray}
&&\frac{1}{4}\frac{d |q|_{4}^{4}}{dt}+3||\nabla q|q  |_{2}^{2}+3||\partial_{\xi}q|q|_{2}^{2}+\beta_{s}| q|_{\xi=1}|_{4}^{4}\nonumber\\
&=&\int_{\mho}Q_{2}q^{3}d\mho
\leq |Q_{2}|_{2}|q^{2}|_{3}^{\frac{3}{2}}\leq   c|Q_{2}|_{2} |q^{2}|_{2}^{ \frac{3}{4}}\|q^{2}\|_{1}^{\frac{3}{4}}\nonumber\\
&=&c|Q_{2}|_{2}|q|_{4}^{\frac{3}{2}}(|q|_{4}^{\frac{3}{2}}+| |\nabla q|q|_{2}^{\frac{3}{4}}+|q\partial_{\xi}q |_{2}^{\frac{3}{4}})\nonumber\\
&\leq& \varepsilon(||\nabla q|q|_{2}^{2}+|q\partial_{\xi}q|_{2}^{2}   )+c|Q_{2}|_{2}^{\frac{8}{5}}|q|_{4}^{\frac{12}{5} }+c|Q_{2}|_{2}|q|_{4}^{3}.
\end{eqnarray}
Since $q^{4}(\theta,\varphi, \xi)=-\int_{\xi}^{1}\partial_{\xi'} q^{4}d\xi'+q^{4}|_{\xi=1},$ by H$\mathrm{\ddot{o}}$lder inequality we get
\begin{eqnarray*}
|q|_{4}^{4}\leq c| q|\partial_{\xi}q||_{2}^{2}+\frac{1}{2}|q|_{4}^{4}+ |q|_{\xi=1}|_{4}^{4},
\end{eqnarray*}
which implies that there exists a positive constant $c_{1}$ such that
\begin{eqnarray*}
\frac{d}{dt}|q|_{4}^{4}+c_{1}|q|_{4}^{4}\leq c|Q_{2}|_{2}^{\frac{8}{5}}|q|_{4}^{\frac{12}{5} }+c|Q_{2}|_{2}|q|_{4}^{3}.
\end{eqnarray*}
That is
\begin{eqnarray*}
\frac{d}{dt}|q|_{4}^{2}+c|q|_{4}^{2}\leq c|Q_{2}|_{2}^{2}.
\end{eqnarray*}
Therefore, we have
\begin{eqnarray}
|q(t)|_{4}^{2}\leq e^{-ct}|q_{0}|_{4}^{2}+c|Q_{2}|_{2}^{2}.
\end{eqnarray}
Combining $(3.36)$ and $(3.37),$ we arrive at
\begin{eqnarray}
\int_{t_{0}}^{t_{0}+1}(||\nabla q(t)|q(t)  |_{2}^{2}+||\partial_{\xi}q(t)|q(t)|_{2}^{2}+\beta_{s}| q|_{\xi=1}(t)|_{4}^{4})dt\leq ce^{-c_{1}t_{0}}|q_{0}|_{4}^{2}+c|Q_{2}|_{2}^{2}.
\end{eqnarray}

\subsection{$L^{4}$ estimate for $T$ }
Taking inner product of $(1.12)$ with $T^{3}$ in $L^{2}(\mho)$ yields,
\begin{eqnarray}
&&\frac{1}{4}\frac{|T|_{4}^{4}}{dt}+3||\nabla T|T |_{2}^{2}+3||\partial_{\xi}T|T|_{2}^{2}+\alpha_{s}|T|_{\xi=1}|_{4}^{4}\nonumber\\
&&=\int_{\mho}\frac{bP}{p}\Big{(}\int_{\xi}^{1} \mathrm{div}  v d\xi'  \Big{)}|T|^{2}Td\mho+\int_{\mho}\frac{abP}{p}\Big{(}\int_{\xi}^{1} \mathrm{div}  v d\xi'  \Big{)}q|T|^{2}Td\mho\nonumber\\
&&-\int_{\mho}\Big{[} \nabla_{v}T+ \Big{(}\int_{\xi}^{1} \mathrm{div}  v d\xi'  \Big{)}\partial_{\xi}T   \Big{]}|T|^{2}Td\mho+\int_{\mho}Q_{1}|T|^{2}Td\mho.
\end{eqnarray}
Using the estimate in $\cite{GH2}$,  we know that
\begin{eqnarray}
&&|\int_{\mho}\frac{bP}{p}\Big{(}\int_{\xi}^{1} \mathrm{div} d\xi'   \Big{)}|T|^{2}T |\nonumber\\
&&\leq c(|\nabla_{e_{\theta}} v|_{2}+ |\nabla_{e_{\varphi}} v|_{2} )|T|_{4}^{2}(|\nabla T|_{2}+|T|_{2} ).
\end{eqnarray}
Taking an similar argument as in $\cite{GH2}$,  we get
\begin{eqnarray*}
&&|\int_{\mho}\frac{abP}{p}q\Big{(}\int_{\xi}^{1} \mathrm{div} v d\xi'   \Big{)}T^{3}  |\\
&&\leq c\int_{0}^{1}\Big{[} \Big{(} \int_{S^{2}}q^{4}dS^{2}  \Big{)}^{\frac{1}{4}}  \Big{(} \int_{S^{2}}T^{12}dS^{2}  \Big{)}^{\frac{1}{4}}\Big{]}d\xi
\sup\limits_{\xi\in [0,1]}|\frac{bP}{p}\int_{\xi}^{1}\mathrm{div} vd\xi'|_{L^{2}(S^{2})}.
\end{eqnarray*}
Therefore,
\begin{eqnarray}
&&|\int_{\mho}\frac{abP}{p}q\Big{(}\int_{\xi}^{1} \mathrm{div} v d\xi'   \Big{)}T^{3}  |\nonumber\\
&\leq& c\int_{0}^{1}|q|_{L^{4}(S^{2}) }|T^{2} |_{L^{6}(S^{2}) }^{\frac{3}{2}}d\xi |\mathrm{div} v|_{2}\nonumber\\
&\leq& c\int_{0}^{1}|q|_{L^{4}(S^{2}) }|T^{2} |_{L^{2}(S^{2}) }^{\frac{1}{2}}|\nabla T^{2} |_{L^{2}(S^{2}) }d\xi |\mathrm{div} v|_{2}\nonumber\\
&\leq &\varepsilon |T|\nabla T| |_{2}^{2}+c(|\nabla_{e_{\theta}}v |_{2}^{2}+ |\nabla_{e_{\varphi}}v |_{2}^{2}    )|T|_{4}^{2}|q|_{4}^{2}.
\end{eqnarray}
To estimate the last term on the right of $(3.39),$ we use interpolation inequality and H$\mathrm{\ddot{o}}$lder inequality to have
\begin{eqnarray}
\int_{\mho}Q_{1}T^{3}d\mho&\leq& |Q_{1}|_{2}|T|_{6}^{3}=|Q_{1}|_{2}|T^{2}|_{3}^{\frac{3}{2}}\nonumber\\
&\leq&c|Q_{1}|_{2}|T^{2}|_{2}^{\frac{3}{4}}(|T|\nabla T|  |_{2}^{\frac{3}{4}}+|T|\partial_{\xi}T||_{2}^{\frac{3}{4}}+|T^{2}|_{2}^{\frac{3}{4}}   )\nonumber\\
&\leq& \varepsilon ( |T|\nabla T|  |_{2}^{2}+ |T|\partial_{\xi}T||_{2}^{2}  )+ c|Q_{1}|_{2} |T|_{4}^{3} + c |Q_{1}|_{2}^{\frac{8}{5}} |T|_{4}^{\frac{12}{5}}.
\end{eqnarray}
By virtue of $(3.39)-(3.42),$ we have
\begin{eqnarray*}
\frac{|T|_{4}^{2}}{2}\frac{d|T|_{4}^{2}}{dt}
&\leq& c(|\nabla_{e_{\theta}} v|_{2}+ |\nabla_{e_{\varphi}} v|_{2} )|T|_{4}^{2}(|\nabla T|_{2}+|T|_{2} )\\
&&+c(|\nabla_{e_{\theta}}v |_{2}^{2}+ |\nabla_{e_{\varphi}}v |_{2}^{2}    )|T|_{4}^{2}|q|_{4}^{2}\\
&&+ c|Q_{1}|_{2} |T|_{4}^{3} + c |Q_{1}|_{2}^{\frac{8}{5}} |T|_{4}^{\frac{12}{5}} ,
\end{eqnarray*}
which implies
\begin{eqnarray}
\frac{d|T|_{4}^{2}}{dt}&\leq& c(|\nabla_{e_{\theta}} v|_{2}+ |\nabla_{e_{\varphi}} v|_{2} )(|\nabla T|_{2}+|T|_{2} )\nonumber\\
&&+c(|\nabla_{e_{\theta}}v |_{2}^{2}+ |\nabla_{e_{\varphi}}v |_{2}^{2}    )|q|_{4}^{2}\nonumber\\
&&+ c|Q_{1}|_{2} |T|_{4} + c |Q_{1}|_{2}^{\frac{8}{5}} |T|_{4}^{\frac{2}{5}} .
\end{eqnarray}
Since
$
|T|_{4}\leq c\|T\|_{1}
$
, by $(3.32), (3.33), (3.37)$ and the uniform Gronwall lemma, we obtain the desired uniform boundedness of   $|T|_{4}$ and it gives the absorbing ball of $T$ in $L^{4}(\mho).$ That is to say there exists a  constant $c$ independent of $t$ such that
\begin{eqnarray}
|T|_{4}(t)\leq c
\end{eqnarray}
for all $t\geq 0.$
Moreover, by $(3.39)$ we have the following uniform bound on the time average
\begin{eqnarray}
&&\int_{t_{0}}^{t_{0}+1}(||\nabla T(t)|T(t) |_{2}^{2}+||\partial_{\xi}T(t)|T(t)|_{2}^{2}+|T(t)|_{\xi=1}|_{4}^{4})dt\nonumber\\
&\leq& c|T(t_{0})|_{4}^{4} +c\int_{t_{0}}^{t_{0}+1}(|\nabla_{e_{\theta}} v(t)|_{2}+ |\nabla_{e_{\varphi}} v(t)|_{2} )|T(t)|_{4}^{2}(|\nabla T(t)|_{2}+|T(t)|_{2} )dt\nonumber\\
&&+c\int_{t_{0}}^{t_{0}+1}(|\nabla_{e_{\theta}}v(t) |_{2}^{2}+ |\nabla_{e_{\varphi}}v(t) |_{2}^{2}    )|T(t)|_{4}^{2}|q(t)|_{4}^{2}dt\nonumber\\
&&+ c\int_{t_{0}}^{t_{0}+1}|Q_{1}|_{2} |T(t)|_{4}^{3}dt + c\int_{t_{0}}^{t_{0}+1} |Q_{1}|_{2}^{\frac{8}{5}} |T(t)|_{4}^{\frac{12}{5}}dt\leq c ,
\end{eqnarray}
where the constant $c$ is independent of $t_{0}.$
\begin{Rem}
In this part, we try to use uniform Gronwall lemma to get the uniform boundedness of $|T|_{4}^{2}$ instead of $|T|_{4}^{4}$. Therefore, from $(3.43)-(3.45)$ we can see that the estimate of $|T|_{3}^{3}$ is not necessary , which simplify the proof of the existence of the absorbing ball in space $L^{4}(\mho)$ in $\cite{GH2}.$ This technique will be used again in the following to get the uniform estimates for $v$ in $L^{4}(\mho)$ with respect to time.
\end{Rem}

\subsection{$L^{4}$ estimate for $v$ }
Taking inner product of equation $(2.22)$ with $|\tilde{v}|^{2}\tilde{v}$ in $L^{2}(\mho)$, we get
\begin{eqnarray}
\frac{1}{4}\frac{d|\tilde{v}|_{4}^{4}}{dt}&+&\int_{\mho}\Big{(}|\nabla_{e_{\theta}}\tilde{v} |^{2}|\tilde{v}|^{2}+|\nabla_{e_{\varphi}}\tilde{v} |^{2} |\tilde{v}|^{2}+\frac{1}{2}|\nabla_{e_{\theta}}|\tilde{v}|^{2} |^{2}+\frac{1}{2}|\nabla_{e_{\varphi}}|\tilde{v}|^{2} |^{2}+|\tilde{v}|^{4}\Big{)}d\mho\nonumber\\
&&+\int_{\mho}\Big{(}|\tilde{v}_{\xi}|^{2}|\tilde{v}|^{2}+\frac{1}{2}|\partial_{\xi}|\tilde{v}|^{2}|^{2}    \Big{)}d\mho
=-\int_{\mho}\Big{[}\nabla_{\tilde{v}}\tilde{v} + \Big{(} \int_{\xi}^{1}\mathrm{div} \tilde{v}d\xi'   \Big{)} \partial_{\xi}\tilde{v}  \Big{]}\cdot|\tilde{v}|^{2}\tilde{v} d\mho\nonumber\\
&&-\int_{\mho}(\nabla_{\bar{v}}\tilde{v})\cdot |\tilde{v}|^{2}\tilde{v}d\mho-\int_{\mho}(\nabla_{\tilde{v}}\bar{v})\cdot |\tilde{v}|^{2}\tilde{v}d\mho\nonumber\\
&&-\int_{\mho}\Big{(}\int_{\xi}^{1}\frac{bP}{p}\mathrm{grad}[(1+aq )T]d\xi'\cdot|\tilde{v}|^{2}\tilde{v}\Big{)}d\mho\nonumber\\
&&+\int_{\mho}\Big{(}\int_{0}^{1}\int_{\xi}^{1}\frac{bP}{p}\mathrm{grad} [(1+aq  )T]d\xi'd\xi\Big{)}\cdot|\tilde{v}|^{2}\tilde{v}d\mho\nonumber\\
&&+\int_{\mho}\overline{(\tilde{v}\mathrm{div} \tilde{v} +\nabla_{\tilde{v}}\tilde{v} )}\cdot |\tilde{v}|^{2}\tilde{v}d\mho-\int_{\mho}(\frac{f}{R_{0}}k\times \tilde{v}   )\cdot |\tilde{v}|^{2}\tilde{v}d\mho,
\end{eqnarray}
where $\tilde{v}_{\xi}=\partial_{\xi}\tilde{v}.$  By Lemma $2.3,$ we have
\begin{eqnarray}
\int_{\mho}\Big{[} \nabla_{\tilde{v}}\tilde{v}+\Big{(} \int_{\xi}^{1}\mathrm{div} \tilde{v} d\xi' \Big{)} \partial_{\xi}\tilde{v}  \Big{]}|\tilde{v}|^{2}\tilde{v}d\mho=0.
\end{eqnarray}
Using Lemma $2.3$ again and $(2.18)$, we obtain
\begin{eqnarray}
\int_{\mho}(\nabla_{\bar{v}}\tilde{v}  )\cdot|\tilde{v}|^{2}\tilde{v}d\mho=\frac{1}{4}\int_{\mho}\nabla_{\bar{v}}|\tilde{v}|^{4}d\mho=-\frac{1}{4}\int_{\mho}|\tilde{v}|^{4}\mathrm{div} \bar{v}d\mho =0.
\end{eqnarray}
By virtue of Lemma $2.2$, we have
\begin{eqnarray*}
0&=&\int_{\mho}\mathrm{div}[(|\tilde{v}|^{2}\tilde{v}\cdot \bar{v}  )\tilde{v} ]d\mho=\int_{\mho}\nabla_{\tilde{v}}(|\tilde{v}|^{2}\tilde{v}\cdot \bar{v} )d\mho+\int_{\mho}|\tilde{v}|^{2}\tilde{v}\cdot\bar{v} \mathrm{div} \tilde{v}d\mho\nonumber\\
&=&\int_{\mho}[|\tilde{v}|^{2}\tilde{v}\cdot\nabla_{\tilde{v}}\bar{v}+\bar{v}\cdot\nabla_{\tilde{v}}(|\tilde{v}|^{2}\tilde{v} )  ]d\mho+\int_{\mho}|\tilde{v}|^{2}\tilde{v}\cdot\bar{v} \mathrm{div} \tilde{v}d\mho.
\end{eqnarray*}
Therefore,
\begin{eqnarray}
\int_{\mho}[|\tilde{v}|^{2}\tilde{v}\cdot\bar{v} \mathrm{div} \tilde{v}+\bar{v}\cdot\nabla_{\tilde{v}}(|\tilde{v}|^{2}\tilde{v} )  ]d\mho=-\int_{\mho}|\tilde{v}|^{2}\tilde{v}\cdot\nabla_{\tilde{v}}\bar{v}d\mho.
\end{eqnarray}
Using integration by parts,  we obtain
\begin{eqnarray}
\int_{\mho}\Big{[} \int_{0}^{1}(\tilde{v} \mathrm{div} \tilde{v} +\nabla_{\tilde{v}}\tilde{v}  )d\xi    \Big{]}\cdot |\tilde{v}|^{2}\tilde{v}d\mho&=&
-\int_{\mho}\Big{(}\int_{0}^{1}\tilde{v}_{\theta}\tilde{v}d\xi   \Big{)}\cdot\nabla_{e_{\theta}}(|\tilde{v}|^{2}\tilde{v})d\mho\nonumber\\
&&-\int_{\mho}\Big{(}\int_{0}^{1}\tilde{v}_{\varphi}\tilde{v}d\xi   \Big{)}\cdot\nabla_{e_{\varphi}}(|\tilde{v}|^{2}\tilde{v} )d\mho.
\end{eqnarray}
Note that the minus on the right hand side of $(3.50)$ was missed in $\cite{GH2}.$  In view of $(3.46)-(3.50)$ combined with $(\frac{f}{R_{0}}\times \tilde{v})\cdot |\tilde{v}|^{2}\tilde{v}=0$, we have
\begin{eqnarray*}
\frac{1}{4}\frac{d|\tilde{v}|_{4}^{4}}{dt}&+&\int_{\mho}\Big{(}|\nabla_{e_{\theta}}\tilde{v} |^{2}|\tilde{v}|^{2}+|\nabla_{e_{\varphi}}\tilde{v} |^{2} |\tilde{v}|^{2}+\frac{1}{2}|\nabla_{e_{\theta}}|\tilde{v}|^{2} |^{2}+\frac{1}{2}|\nabla_{e_{\varphi}}|\tilde{v}|^{2} |^{2}+|\tilde{v}|^{4}\Big{)}d\mho\nonumber\\
&&+\int_{\mho}\Big{(}|\tilde{v}_{\xi}|^{2}|\tilde{v}|^{2}+\frac{1}{2}|\partial_{\xi}|\tilde{v}|^{2}|^{2}    \Big{)}d\mho\nonumber\\
&&=\int_{\mho}[|\tilde{v}|^{2}\tilde{v}\cdot\bar{v} \mathrm{div} \tilde{v}+\bar{v}\cdot\nabla_{\tilde{v}}(|\tilde{v}|^{2}\tilde{v} )  ]d\mho\nonumber\\
&&+\int_{\mho}\Big{(}\int_{0}^{1}\tilde{v}_{\theta}\tilde{v}d\xi   \Big{)}\cdot\nabla_{e_{\theta}}(|\tilde{v}|^{2}\tilde{v})d\mho
+\int_{\mho}\Big{(}\int_{0}^{1}\tilde{v}_{\varphi}\tilde{v}d\xi   \Big{)}\cdot\nabla_{e_{\varphi}}(|\tilde{v}|^{2}\tilde{v} )d\mho\nonumber\\
&&+\int_{\mho}\Big{(}\int_{\xi}^{1}\frac{bP}{p}[(1+aq )T]d\xi'\cdot\mathrm{div}(|\tilde{v}|^{2}\tilde{v})\Big{)}d\mho\nonumber\\
&&-\int_{\mho}\Big{(}\int_{0}^{1}\int_{\xi}^{1}\frac{bP}{p}[(1+aq )T]d\xi'd\xi\cdot\mathrm{div}(|\tilde{v}|^{2}\tilde{v})\Big{)}d\mho.
\end{eqnarray*}
Then by H$\mathrm{\ddot{o}}$lder inequality and the argument above, we obtain
\begin{eqnarray}
\frac{1}{4}\frac{d|\tilde{v}|_{4}^{4}}{dt}&+&\int_{\mho}\Big{(}|\nabla_{e_{\theta}}\tilde{v} |^{2}|\tilde{v}|^{2}+|\nabla_{e_{\varphi}}\tilde{v} |^{2} |\tilde{v}|^{2}+\frac{1}{2}|\nabla_{e_{\theta}}|\tilde{v}|^{2} |^{2}+\frac{1}{2}|\nabla_{e_{\varphi}}|\tilde{v}|^{2} |^{2}+|\tilde{v}|^{4}\Big{)}d\mho\nonumber\\
&&+\int_{\mho}\Big{(}|\tilde{v}_{\xi}|^{2}|\tilde{v}|^{2}+\frac{1}{2}|\partial_{\xi}|\tilde{v}|^{2}|^{2}    \Big{)}d\mho\nonumber\\
&\leq&c\int_{S^{2}}|\bar{v}|\Big{(}\int_{0}^{1}|\tilde{v}|^{3}(|\nabla_{e_{\theta}}\tilde{v} |^{2}+ |\nabla_{e_{\varphi}}\tilde{v} |^{2}  )^{\frac{1}{2}}d\xi\Big{)} dS^{2}\nonumber\\
&&+c\int_{S^{2}}\Big{(}\int_{0}^{1}|\tilde{v}|^{2}d\xi\Big{)}\Big{(}\int_{0}^{1}|\tilde{v}|^{2}(|\nabla_{e_{\theta}}\tilde{v} |^{2}+ |\nabla_{e_{\varphi}}\tilde{v} |^{2}  )^{\frac{1}{2}}d\xi\Big{)} dS^{2}\nonumber\\
&&+c\int_{S^{2}}\Big{(} (|\bar{T}|+|\overline{qT}|   )\int_{0}^{1}|\tilde{v}|^{2}(|\nabla_{e_{\theta}}\tilde{v} |^{2}+ |\nabla_{e_{\varphi}}\tilde{v} |^{2}  )^{\frac{1}{2}}   d\xi  \Big{)}dS^{2}\nonumber\\
&=:&I_{1}+I_{2}+I_{3}.
\end{eqnarray}
In the following, we estimate $I_{i},i=1,2,3,$ separately. By interpolation inequality , H$\mathrm{\ddot{o}}$lder inequality and Minkowski inequality, we have
\begin{eqnarray*}
I_{1}&\leq&\int_{S^{2}}\Big{[}|\bar{v}| \Big{(}\int_{0}^{1}|\tilde{v}|^{4}d\xi\Big{)}^{\frac{1}{2}} \Big{(} \int_{0}^{1}|\tilde{v}|^{2}(|\nabla_{e_{\theta}}\tilde{v}|^{2}+  |\nabla_{e_{\varphi}}\tilde{v}|^{2} )d\xi   \Big{)}^{\frac{1}{2}}  \Big{]}dS^{2}\nonumber\\
&\leq&|\bar{v}|_{4}\Big{(}\int_{S^{2}}(\int_{0}^{1}|\tilde{v}|^{4}d\xi   )^{2} dS^{2}  \Big{)}^{\frac{1}{4}}
\Big{(} \int_{\mho}|\tilde{v}|^{2}(|\nabla_{e_{\theta}}\tilde{v}|^{2}+ |\nabla_{e_{\varphi}}\tilde{v}|^{2}) d\mho  \Big{)}^{\frac{1}{2}}\nonumber\\
&\leq&|\bar{v}|_{4}\Big{(} \int_{0}^{1}( \int_{S^{2}}|\tilde{v}|^{2\times4}dS^{2} )^{\frac{1}{2}}   d\xi \Big{)}^{\frac{1}{2}}\Big{(} \int_{\mho}|\tilde{v}|^{2}(|\nabla_{e_{\theta}}\tilde{v}|^{2}+ |\nabla_{e_{\varphi}}\tilde{v}|^{2})d\mho   \Big{)}^{\frac{1}{2}}\nonumber\\
&\leq&\varepsilon \int_{\mho}|\tilde{v}|^{2}(|\nabla_{e_{\theta}}\tilde{v}|^{2}+ |\nabla_{e_{\varphi}}\tilde{v}|^{2})d\mho +c|\bar{v}|_{4}^{2}\int_{0}^{1}||\tilde{v}|^{2}|_{2}(||\tilde{v}|^{2}  |_{2}+|\nabla |\tilde{v}|^{2} |_{2}  )d\xi\nonumber\\
&\leq&\varepsilon \int_{\mho}|\tilde{v}|^{2}(|\nabla_{e_{\theta}}\tilde{v}|^{2}+ |\nabla_{e_{\varphi}}\tilde{v}|^{2})d\mho +c|\bar{v}|_{4}^{2}
||\tilde{v}|^{2}|_{2}(||\tilde{v}|^{2}  |_{2}+|\nabla |\tilde{v}|^{2} |_{2}  )\nonumber\\
&\leq&\varepsilon \int_{\mho}|\tilde{v}|^{2}(|\nabla_{e_{\theta}}\tilde{v}|^{2}+|\nabla_{e_{\varphi}}\tilde{v}|^{2})d\mho +\varepsilon |\nabla |\tilde{v}| ^{2} |_{2}^{2}+c(\|\bar{v}\|_{1}^{2}+|\bar{v}|_{2}^{2}\|\bar{v}\|_{1}^{2} )|\tilde{v}|_{4}^{4}.
\end{eqnarray*}
Similarly,
\begin{eqnarray*}
I_{2}&\leq& c\int_{S^{2}}( \int_{0}^{1}|\tilde{v}|^{2}d\xi )( \int_{0}^{1}|\tilde{v}|^{2}d\xi )^{\frac{1}{2}}\Big{(} \int_{0}^{1}|\tilde{v}|^{2}(|\nabla_{e_{\theta}}\tilde{v}|^{2}+|\nabla_{e_{\varphi}}\tilde{v} |^{2}  )d\xi  \Big{)}^{\frac{1}{2}}dS^{2}\nonumber\\
&\leq&c\Big{(}\int_{\mho}|\tilde{v}|^{2}(|\nabla_{e_{\theta}}\tilde{v}|^{2}+|\nabla_{e_{\varphi}}\tilde{v} |^{2})d\mho    \Big{)}^{\frac{1}{2}}
\Big{(}\int_{S^{2}} (\int_{0}^{1}|\tilde{v}|^{2}d\xi   )^{3}dS^{2}    \Big{)}^{\frac{1}{2}}\nonumber\\
&\leq&\varepsilon\int_{\mho}|\tilde{v}|^{2}(|\nabla_{e_{\theta}}\tilde{v}|^{2}+|\nabla_{e_{\varphi}}\tilde{v} |^{2})d\mho
+c\Big{(}\int_{0}^{1}(\int_{S^{2}}|\tilde{v}|^{6}dS^{2}  )^{\frac{1}{3}}d\xi    \Big{)}^{3}\nonumber\\
&\leq&\varepsilon\int_{\mho}|\tilde{v}|^{2}(|\nabla_{e_{\theta}}\tilde{v}|^{2}+|\nabla_{e_{\varphi}}\tilde{v} |^{2})d\mho +
c\Big{(}\int_{0}^{1}|\tilde{v}|_{4}^{\frac{4}{3}}(|\tilde{v}|_{2}^{\frac{2}{3}}+|\nabla_{e_{\theta}} \tilde{v}|_{2}^{\frac{2}{3}}+ |\nabla_{e_{\varphi}} \tilde{v}|_{2}^{\frac{2}{3}}    )d\xi    \Big{)}^{3}\nonumber\\
&\leq&\varepsilon\int_{\mho}|\tilde{v}|^{2}(|\nabla_{e_{\theta}}\tilde{v}|^{2}+|\nabla_{e_{\varphi}}\tilde{v} |^{2})d\mho+c|\tilde{v}|_{4}^{4}|\tilde{v}|_{2}^{2}+c|\tilde{v}|_{4}^{4}(|\nabla_{e_{\theta}}\tilde{v}|_{2}^{2}+ |\nabla_{e_{\varphi}}\tilde{v}|_{2}^{2} ).
\end{eqnarray*}
Analogously, we have
\begin{eqnarray*}
I_{3}&\leq& \int_{S^{2}}(|\overline{T}|+|\overline{qT}|   )\Big{(}\int_{0}^{1}|\tilde{v}|^{2}d\xi  \Big{)}^{\frac{1}{2}}\Big{(} \int_{0}^{1}|\tilde{v}|^{2}(|\nabla_{e_{\theta}} \tilde{v}|^{2}+ |\nabla_{e_{\varphi}} \tilde{v}|^{2}  ) d\xi \Big{)}^{\frac{1}{2}}dS^{2}\nonumber\\
&\leq& \Big{(}\int_{\mho}|\tilde{v}|^{2}(|\nabla_{e_{\theta}} \tilde{v}|^{2}+ |\nabla_{e_{\varphi}} \tilde{v}|^{2}  )d\mho\Big{)}^{\frac{1}{2}}
\Big{(}\int_{S^{2}}(\int_{0}^{1}|\tilde{v}|^{2}d\xi   )^{2}  dS^{2}\Big{)}^{\frac{1}{4}}(|\overline{T}|_{4}+|\overline{qT}|_{4}  )\nonumber\\
&\leq& \varepsilon \int_{\mho}|\tilde{v}|^{2}(|\nabla_{e_{\theta}} \tilde{v}|^{2}+ |\nabla_{e_{\varphi}} \tilde{v}|^{2}  )d\mho
+c\int_{0}^{1}\Big{(} \int_{S^{2}}|\tilde{v}|^{4}dS^{2}   \Big{)}^{\frac{1}{2}}d\xi( |\overline{T}|_{4}^{2}+|\overline{qT}|_{4}^{2}  )\nonumber\\
&\leq& \varepsilon \int_{\mho}|\tilde{v}|^{2}(|\nabla_{e_{\theta}} \tilde{v}|^{2}+ c|\nabla_{e_{\varphi}} \tilde{v}|^{2}  )d\mho+|\tilde{v}|_{4}^{2}
\Big{[}|\overline{T}|_{4}^{2}+\Big{(}\int_{S^{2}}(\int_{0}^{1}qT d\xi )^{4}dS^{2}\Big{)}^{\frac{1}{2}} \Big{]} \nonumber\\
&\leq& \varepsilon \int_{\mho}|\tilde{v}|^{2}(|\nabla_{e_{\theta}} \tilde{v}|^{2}+ c|\nabla_{e_{\varphi}} \tilde{v}|^{2}  )d\mho
+ c |\tilde{v}|_{4}^{2}|T|_{4}^{2} +c|\tilde{v}|_{4}^{2}\Big{(}\int_{0}^{1}|q|_{8}|T|_{8}d\xi\Big{)}^{2}\nonumber\\
&\leq& \varepsilon \int_{\mho}|\tilde{v}|^{2}(|\nabla_{e_{\theta}} \tilde{v}|^{2}+ c|\nabla_{e_{\varphi}} \tilde{v}|^{2}  )d\mho
+c |\tilde{v}|_{4}^{2}|T|_{4}^{2}    \nonumber\\
&&+c|\tilde{v}|_{4}^{2}\Big{(}\int_{0}^{1}|q|_{4}^{\frac{1}{2}}(|q|_{2}^{\frac{1}{2}}+|\nabla q|_{2}^{\frac{1}{2}})|T|_{4}^{\frac{1}{2}}( |T|_{2}^{\frac{1}{2}}+  |\nabla T|_{2}^{\frac{1}{2}} )d\xi \Big{)}^{2}\nonumber\\
&\leq& \varepsilon \int_{\mho}|\tilde{v}|^{2}(|\nabla_{e_{\theta}} \tilde{v}|^{2}+ c|\nabla_{e_{\varphi}} \tilde{v}|^{2}  )d\mho
+c|\tilde{v}|_{4}^{2}|T|_{4}^{2}\nonumber\\
&&+c|\tilde{v}|_{4}^{2}(|q|_{4}^{2}+|q|_{4}|\nabla q|_{2}  )(         |T|_{4}^{2}  +|T|_{4}  |\nabla T|_{2}).
\end{eqnarray*}
By the estimates of $I_{1}-I_{3},$ we have
\begin{eqnarray*}
\frac{d|\tilde{v}|_{4}^{2}}{dt}\leq c(\|v \|_{1}^{2}+|v|_{2}^{2}|\|v\|_{1}^{2} )|\tilde{v}|_{4}^{2}+c|T|_{4}^{2}+c(|q|_{4}^{2}+|q|_{4}|\nabla q|_{2}  )(         |T|_{4}^{2}  +|T|_{4}  |\nabla T|_{2}).
\end{eqnarray*}
Therefore, by $(3.27), (3.32), (3.33), (3.37), (3.44)$ and uniformly Gronwall Lemma, we have
\begin{eqnarray}
|\tilde{v}(t)|_{4}\leq c,
\end{eqnarray}
where $t\geq 0$ and $c$ is independent of $t$. Furthermore, by $(3.51)$ and estimates of $I_{1}-I_{3},$ we get
\begin{eqnarray}
&&\int_{t_{0}}^{t_{0}+1} ||\nabla_{e_{\theta}}\tilde{v}(t) ||\tilde{v}(t)||_{2}^{2}+||\nabla_{e_{\varphi}}\tilde{v}(t) | |\tilde{v}(t)||_{2}^{2}dt\nonumber\\
&\leq & |\tilde{v}(t_{0})|_{4}^{4}+c \int_{t_{0}}^{t_{0}+1}[(\|v(t) \|_{1}^{2}+|v(t)|_{2}^{2}|\|v(t)\|_{1}^{2} )|\tilde{v}(t)|_{4}^{2}+|T(t)|_{4}^{2}]dt\nonumber\\
&&+c\int_{t_{0}}^{t_{0}+1} (|q(t)|_{4}^{2}+|q(t)|_{4}|\nabla q(t)|_{2}  )(         |T(t)|_{4}^{2}  +|T(t)|_{4}  |\nabla T(t)|_{2})dt\leq c,
\end{eqnarray}
where $c$ is independent of $t_{0}.$
\subsection{$H^{1}$ estimates of $v,T,q$. }
From $\cite{GH2},$ we have
\begin{eqnarray}
\frac{d\|\bar{v}\|_{1}}{dt}+|\Delta \bar{v}|_{2}^{2}\leq c(\|\bar{v}\|_{1}^{2}+|\bar{v}|_{2}^{2}\|\bar{v}\|_{1}^{2}  )\|\bar{v}\|_{1}^{2}+c||\tilde{v} |  |\nabla_{e_{\theta}}\tilde{v}||_{2}^{2}+c||\tilde{v} |  |\nabla_{e_{\varphi}}\tilde{v}||_{2}^{2}.
\end{eqnarray}
By the uniform Gronwall lemma and $(3.32)-(3.33)$, we obtain
\begin{eqnarray}
\|\bar{v}(t)\|_{1}\leq c,
\end{eqnarray}
where $t\geq 0$ and $c$ is independent of $t.$ By Sobolev inequality, we have that for all $t\geq 0$
\begin{eqnarray*}
|\bar{v}(t)|_{L^{4}(S^{2})}
\leq c|\bar{v}(t)|_{L^{2}(S^{2})}+c|\nabla \bar{v}(t)|_{L^{2}(S^{2})}\leq c,
\end{eqnarray*}
which implies
\begin{eqnarray}
|v(t)|_{4}  \leq |\tilde{v}(t)|_{4} + |\bar{v}(t)|_{4}\leq c.
\end{eqnarray}
From $\cite{GH2},$ we have
\begin{eqnarray}
&&\frac{d|v_{\xi}|_{2}^{2}}{dt}+|\nabla_{e_{\theta}}v_{\xi} |_{2}^{2}+|\nabla_{e_{\varphi}}v_{\xi} |_{2}^{2}+|v_{\xi}|_{2}^{2}+|v_{\xi\xi}|_{2}^{2}\nonumber\\
&\leq&c(\|\bar{v}\|_{1}^{8}+|\tilde{v}|_{4}^{8}  )|v_{\xi}|_{2}^{2}+c|T|_{2}^{2}+c|q|_{4}^{4}+c|T|_{4}^{4},
\end{eqnarray}
where $ \partial_{\xi\xi} v= v_{\xi\xi}.$ By the uniform Gronwall lemma, we obtain
\begin{eqnarray}
|v_{\xi}(t)|_{2}\leq c,
\end{eqnarray}
where $t\geq 0$ and $c$ is independent of $t.$ Therefore, by $(3.56)-(3.58)$ and uniformly Gronwall inequality we have
\begin{eqnarray}
\int_{t_{0}}^{t_{0}+1} ( |\nabla_{e_{\theta}}v_{\xi} |_{2}^{2}+|\nabla_{e_{\varphi}}v_{\xi} |_{2}^{2}+|v_{\xi\xi}|_{2}^{2} )dt\leq c,
\end{eqnarray}
where $t_{0}\geq 0$ and $c$ is independent of $t_{0}.$
Taking inner product with $-\Delta v$ in $L^{2}(\mho),$ we obtain
\begin{eqnarray}
&&\frac{1}{2}\frac{d(|\nabla_{e_{\theta}} v|_{2}^{2}+ |\nabla_{e_{\varphi}} v|_{2}^{2}+|v|_{2}^{2}   )}{dt}+|\Delta v|_{2}^{2}+|\nabla_{e_{\theta}}v_{\xi} |_{2}^{2}+|\nabla_{e_{\varphi}}v_{\xi} |_{2}^{2}+|v_{\xi}|_{2}^{2}\nonumber\\
&=&\int_{\mho}\nabla_{v}v\cdot\Delta vd\mho+\int_{\mho}\Big{(}\int_{\xi}^{1}\mathrm{div} vd\xi'\Big{)}\partial_{\xi}v \cdot \Delta vd\mho\nonumber\\
&&+\int_{\mho}\Big{(}\int_{\xi}^{1}\frac{bP}{p}\mathrm{grad} (1+aq)Td\xi'\Big{)}\Delta vd\mho\nonumber\\
&&+\int_{\mho} \frac{f}{R_{0}}v^{\bot}\cdot\Delta vd\mho+\int_{\mho} \mathrm{grad}\Phi_{s}\cdot\Delta vd\mho\nonumber\\
&=& J_{1}+J_{2}+J_{3}+J_{4}+J_{5}.
\end{eqnarray}
By H$\mathrm{\ddot{o}}$lder inequality, interpolation inequality and Young's inequality, we have
\begin{eqnarray*}
J_{1}&\leq& |\Delta v|_{2}|v|_{4}(|\nabla_{e_{\theta}} v|_{4}+ |\nabla_{e_{\varphi}} v|_{4})\nonumber\\
&\leq& c|\Delta v|_{2}|v|_{4}(|\nabla_{e_{\theta}} v|_{2}^{\frac{1}{4}}+|\nabla_{e_{\varphi}} v |_{2}^{\frac{1}{4}}     )
(|\Delta v |_{2}^{\frac{3}{4}}+ |\nabla_{e_{\theta}} v_{\xi} |_{2}^{\frac{3}{4}}+ |\nabla_{e_{\varphi}} v_{\xi} |_{2}^{\frac{3}{4}} +|\nabla_{e_{\varphi}} v|_{2}^{ \frac{3}{4}}+ |\nabla_{e_{\theta}} v|_{2}^{ \frac{3}{4}})\\
&\leq &\varepsilon |\Delta v|_{2}^{2}+\varepsilon |\nabla_{e_{\theta}} v_{\xi}|_{2}^{2}+\varepsilon |\nabla_{e_{\varphi}} v_{\xi}|_{2}^{2} +c(|\nabla_{e_{\theta}} v|_{2}^{2}+|\nabla_{e_{\varphi}} v|_{2}^{2}  ) (|v|_{4}^{8}+|v|_{4}^{2})\\
&\leq &\varepsilon |\Delta v|_{2}^{2}+\varepsilon |\nabla_{e_{\theta}} v_{\xi}|_{2}^{2}+\varepsilon |\nabla_{e_{\varphi}} v_{\xi}|_{2}^{2} +c(|\nabla_{e_{\theta}} v|_{2}^{2}+|\nabla_{e_{\varphi}} v|_{2}^{2}  ).
\end{eqnarray*}
To estimate $J_{2},$ we have
\begin{eqnarray}
J_{2}&\leq& \int_{S^{2}}\Big{(}\int_{0}^{1}|\mathrm{div} v|d\xi \Big{)} \Big{(}\int_{0}^{1}|\partial_{\xi}v|  |\Delta v|d\xi \Big{)} dS^{2}\nonumber\\
&\leq&\int_{S^{2}}\Big{(}\int_{0}^{1}|\mathrm{div} v|d\xi \Big{)} \Big{(}\int_{0}^{1}|\partial_{\xi}v|^{2}d\xi \Big{)}^{\frac{1}{2}}
\Big{(}\int_{0}^{1}|\Delta v|^{2}d\xi \Big{)}^{\frac{1}{2}}dS^{2}\nonumber\\
&\leq&c|\Delta v|_{2}\Big{(}\int_{S^{2}}( \int_{0}^{1}|\partial_{\xi}v|^{2}d\xi  )^{2} dS^{2}\Big{)}^{\frac{1}{4}}\Big{(}\int_{S^{2}}( \int_{0}^{1}|\mathrm{div} v|d\xi  )^{4}dS^{2} \Big{)}^{\frac{1}{4}}\nonumber\\
&\leq&c|\Delta v|_{2}\Big{(} \int_{0}^{1}( \int_{S^{2}} |\partial_{\xi} v |^{4}dS^{2}  )^{\frac{1}{2}} d\xi        \Big{)}^{\frac{1}{2}}
\Big{(} \int_{0}^{1}( \int_{S^{2}} |\mathrm{div} v |^{4}dS^{2}  )^{\frac{1}{4}} d\xi        \Big{)}\nonumber\\
&\leq& c|\Delta v|_{2}\Big{(} \int_{0}^{1} | v_{\xi}|_{2}(|v_{\xi}|_{2}+| \nabla_{e_{\theta}}  v_{\xi}|_{2}
    +|\nabla_{e_{\varphi}} v_{\xi} |_{2})d\xi        \Big{)}^{\frac{1}{2}}\cdot\nonumber\\
&&\Big{(} \int_{0}^{1} |\mathrm{div} v|_{2}^{\frac{1}{2}}(|\mathrm{div} v|_{2}^{\frac{1}{2}}+| \Delta v|_{2}^{\frac{1}{2}})d\xi        \Big{)}\nonumber\\
&\leq&c |\Delta v|_{2}|v_{\xi}|_{2}^{\frac{1}{2}}(|v_{\xi}|_{2}^{\frac{1}{2}}+|\nabla_{e_{\theta}} v_{\xi}|_{2}^{\frac{1}{2}}+|\nabla_{e_{\varphi}}v_{\xi}|_{2}^{\frac{1}{2}}   )\nonumber\\
&&\cdot |\mathrm{div} v  |_{2}^{ \frac{1}{2}}( |\mathrm{div} v|_{2}^{\frac{1}{2}}+|\Delta v|_{2}^{\frac{1}{2} }  )\nonumber\\
&\leq& \varepsilon |\Delta v|_{2}^{2}+c(|\nabla_{e_{\theta}}v |_{2}^{2}+ |\nabla_{e_{\varphi}}v |_{2}^{2}   )\nonumber\\
&&\cdot(1+ |v_{\xi}|_{2}^{4}+|v_{\xi}|_{2}^{2}|\nabla_{e_{\theta}}v_{\xi}|_{2}^{2}+|v_{\xi}|_{2}^{2}|\nabla_{e_{\varphi}}v_{\xi}|_{2}^{2}  ).
\end{eqnarray}
Concerning estimate of $J_{3},$ there are some typos in $(5.80)$ of $\cite{GH2}$. In the first inequality of $(5.80)$, the integral region should be $S^{2}$ instead of $\Omega.$ For reader's convenience, we estimate $J_{3}$ again. Using H$\mathrm{\ddot{o}}$lder inequality, Minkowski inequality and interpolation inequality, we have
\begin{eqnarray}
J_{3}&\leq& c|\nabla T|_{2}|\Delta v|_{2}+c\int_{S^{2}}(\int_{0}^{1}|\mathrm{grad} q|^{2}d\xi )^{\frac{1}{2}}(\int_{0}^{1}T^{2}d\xi )^{\frac{1}{2}}\int_{0}^{1}|\Delta v|d\xi dS^{2}\nonumber\\
&&+c\int_{S^{2}}(\int_{0}^{1}|\mathrm{grad} T|^{2}d\xi )^{\frac{1}{2}}(\int_{0}^{1}q^{2}d\xi )^{\frac{1}{2}}(\int_{0}^{1}|\Delta v|d\xi) dS^{2}\nonumber\\
&\leq& c|\nabla T|_{2}|\Delta v|_{2}+c|\Delta v|_{2}|T|_{4}\Big{(}\int_{0}^{1}(\int_{S^{2}}|\mathrm{grad} q|^{4} dS^{2} )^{\frac{1}{2}}d\xi    \Big{)}^{\frac{1}{2}}\nonumber\\
&&+c|\Delta v|_{2}|q|_{4}\Big{(}\int_{0}^{1}(\int_{S^{2}}|\mathrm{grad} T|^{4} dS^{2} )^{\frac{1}{2}}d\xi    \Big{)}^{\frac{1}{2}}\nonumber\\
&\leq& c|\nabla T|_{2}|\Delta v|_{2}+c|\Delta v|_{2}|T|_{4}|\nabla q|_{2}^{\frac{1}{2}}(|\nabla q|_{2}^{\frac{1}{2}}   +|\Delta q|_{2}^{\frac{1}{2}} )\nonumber\\
&&+c|\Delta v|_{2}|q|_{4}|\nabla T|_{2}^{\frac{1}{2}}(|\nabla T|_{2}^{\frac{1}{2}}   +|\Delta T|_{2}^{\frac{1}{2}} )\nonumber\\
&\leq& \varepsilon |\Delta v|_{2}^{2}+ \varepsilon |\Delta T|_{2}^{2}+\varepsilon |\Delta q|_{2}^{2}+c|\nabla T|_{2}^{2}+c|q|_{4}^{2}|\nabla T|_{2}^{2}\nonumber\\
&&+c|T|_{4}^{2}|\nabla q|_{2}^{2}+c|q|_{4}^{4}|\nabla T|_{2}^{2}+c|T|_{4}^{4}|\nabla q|_{2}^{2}.
\end{eqnarray}
By H$\mathrm{\ddot{o}}$lder inequality, we have
\begin{eqnarray*}
J_{4}\leq \varepsilon |\Delta v|_{2}^{2}+c|v|_{2}^{2}.
\end{eqnarray*}
By Lemma $2.1,$ we infer that $J_{5}=0.$
From $(3.60)$ and estimates of $J_{1}-J_{5},$ we obtain
\begin{eqnarray}
&&\frac{d(|\nabla_{e_{\theta}} v|_{2}^{2}+ |\nabla_{e_{\varphi}} v|_{2}^{2}+|v|_{2}^{2}   )}{dt}\nonumber\\
&\leq&c(| \nabla_{e_{\theta}}v |_{2}^{2}+ | \nabla_{e_{\varphi}}v |_{2}^{2}  )(1+|v_{\xi}|_{2}^{4}+|v_{\xi}|_{2}^{2}|\nabla_{e_{\theta}} v_{\xi}|_{2}^{2}+
|v_{\xi}|_{2}^{2}|\nabla_{e_{\varphi}} v_{\xi}|_{2}^{2})\nonumber\\
&&+c+c|\nabla T|_{2}^{2}+c|\nabla q|_{2}^{2}.
\end{eqnarray}
In view of $( 3.27), (3.33), (3.58), (3.59)$ and uniform Gronwall lemma, we get that for all $t\geq 0$
\begin{eqnarray}
|\nabla_{e_{\theta}} v(t)|_{2}^{2}+ |\nabla_{e_{\varphi}} v(t)|_{2}^{2}\leq c,
\end{eqnarray}
where $c$ is a positive constant which is independent of $t.$ Furthermore, by virtue of $(3.60)$ we infer that for all $t_{0}\geq 0$
\begin{eqnarray}
\int_{t_{0}}^{t_{0}+1}|\Delta v(t)|_{2}^{2}+|\nabla_{e_{\theta}}v_{\xi}(t) |_{2}^{2}+|\nabla_{e_{\varphi}}v_{\xi} (t)|_{2}^{2}dt\leq c,
\end{eqnarray}
where $c$ is independent of $t_{0}$.
\begin{Rem}
To get the estimates of $T$ in $V_{2}$ space, we should first take derivative with respect to $\xi$ in the temperature equation and then  estimate $T_{\xi}$ in $L^{2}(\mho).$ If we take inner product of equation $(2.12)$ with $-\partial_{\xi\xi}T$ in $L^{2}(\mho),$ it is difficult to obtain the energy estimates of $T_{\xi}$ because of its higher nonlinear structure than the oceanic primitive equation.   Therefore, one can not even prove the global existence of the strong solution in this way.
\end{Rem}
By integration by parts,
\begin{eqnarray*}
\int_{\mho}Q_{1\xi}T_{\xi}d\mho&=&\int_{\mho}[\partial_{\xi}(Q_{1}T_{\xi})-Q_{1}T_{\xi\xi}]d\mho\\
&=&-\alpha_{s}\int_{S^{2}}Q_{1}|_{\xi=1}T|_{\xi=1}dS^{2}-\int_{\mho}Q_{1}T_{\xi\xi}d\mho\\
&\leq&\varepsilon |T_{\xi\xi}|_{2}^{2}+ c|Q_{1}|_{\xi=1}|_{2}^{2}+c|T|_{\xi=1}|_{2}^{2}+|Q_{1}|_{2}^{2}.
\end{eqnarray*}
In view of estimates in $\cite{GH2}$ and the above argument, we have
\begin{eqnarray}
&&\frac{1}{2}\frac{d(|T_{\xi}|_{2}^{2}+\alpha_{s}|T|_{\xi=1}|_{2}^{2} )}{dt}+ |\nabla T_{\xi}|_{2}^{2}+|T_{\xi\xi}|_{2}^{2}+\alpha_{s}|\nabla T|_{\xi=1}|_{2}^{2}\nonumber\\
&\leq& \varepsilon (|T_{\xi\xi}|_{2}^{2}+ |\nabla T_{\xi}|_{2}^{2}+|q_{\xi\xi}|_{2}^{2}+ |\nabla q_{\xi}|_{2}^{2} )+c|T_{\xi}|_{2}^{2}+c|q_{\xi}|_{2}^{2}+c\|v_{\xi}\|_{1}^{2}+c\|v\|_{1}^{2}\nonumber\\
&&+c|\nabla q|_{2}^{2}+|T|_{\xi=1}|_{4}^{4}+c|T|_{\xi=1}|_{2}^{2}+c|Q_{1}|_{\xi=1}||_{2}^{2}+c|Q_{1}|_{2}^{2}+c.
\end{eqnarray}
Similarly, by the estimates in $\cite{GH2}$, we also have
\begin{eqnarray}
&&\frac{1}{2}\frac{d(|q_{\xi}|_{2}^{2}+\beta_{s}|q|_{\xi=1}|_{2}^{2} )}{dt}+ |\nabla q_{\xi}|_{2}^{2}+|q_{\xi\xi}|_{2}^{2}+\beta_{s}|\nabla q|_{\xi=1}|_{2}^{2}\nonumber\\
&\leq& \varepsilon (|q_{\xi\xi}|_{2}^{2}+ |\nabla q_{\xi}|_{2}^{2} )+c|q_{\xi}|_{2}^{2}+c\|v\|_{1}^{2}+c\|v_{\xi}\|_{1}^{2}\nonumber\\
&&+c|q|_{\xi=1}|_{4}^{4}+c|q|_{\xi=1}|_{2}^{2}++c|Q_{2}|_{\xi=1}||_{2}^{2}+c|Q_{2}|_{2}^{2}.
\end{eqnarray}
Combining $(3.66)$ and $(3.67)$ yields,
\begin{eqnarray}
&&\frac{d(|T_{\xi}|_{2}^{2}+|q_{\xi}|_{2}^{2}+\alpha_{s}|T|_{\xi=1}|_{2}^{2}+\beta_{s}|q|_{\xi=1}|_{2}^{2} )}{dt}+ |\nabla T_{\xi}|_{2}^{2}+|T_{\xi\xi}|_{2}^{2}\nonumber\\
&&+ |\nabla q_{\xi}|_{2}^{2}+|q_{\xi\xi}|_{2}^{2}+\alpha_{s}|\nabla T|_{\xi=1}|_{2}^{2}+\beta_{s}|\nabla q|_{\xi=1}|_{2}^{2}\nonumber\\
&&\leq c+c(|T_{\xi}|_{2}^{2}+ |q_{\xi}|_{2}^{2}   )+c\|v_{\xi}\|_{1}^{2}+c\|v\|_{1}^{2}+c\|q\|_{1}^{2}+c|T|_{\xi=1}|_{4}^{4}+c|T|_{\xi=1}|_{2}^{2}\nonumber\\
&&+c|q|_{\xi=1}|_{4}^{4}+c|q|_{\xi=1}|_{2}^{2}+c(|Q_{1}|_{\xi=1}|_{2}^{2}+|Q_{2}|_{\xi=1}|_{2}^{2} )+c(|Q_{1}|_{2}^{2}+|Q_{2}|_{2}^{2} ).
\end{eqnarray}
Then in view of uniform Gronwall lemma, $(3.27), (3.33), (3.38), (3.45)$ and $(3.59),$ we obtain for arbitrary $t_{0}\geq 0$
\begin{eqnarray}
|T_{\xi}(t_{0})|_{2}^{2}+|q_{\xi}(t_{0})|_{2}^{2}+ \int_{t_{0}}^{t_{0}+1}( \|T_{\xi}(t)\|_{1}^{2}+\|q_{\xi}(t)\|_{1}^{2}+|\nabla T|_{\xi=1}(t)|_{2}^{2}+|\nabla q|_{\xi}(t)|_{2}^{2}  )dt\leq c,
\end{eqnarray}
where $c$ is independent of $t_{0}.$ By taking inner product of equation $(1.12)$ with $-\Delta T,$ in $L^{2}(\mho),$ we reach
\begin{eqnarray}
&&\frac{1}{2}\frac{d |\nabla T|_{2}^{2} }{dt}+|\Delta T|_{2}^{2}
+|\nabla T_{z}|_{2}^{2}+\alpha_{s}|\nabla T|_{\xi=1} |_{L^{2}(S^{2})}^{2}\nonumber\\
&=&\int_{\mho}\nabla_{v}T \Delta T d\mho+\int_{\mho}\int_{\xi}^{1} \mathrm{div} vd\xi' T_{\xi} \Delta T d\mho\nonumber\\
&&-\int_{\mho}\frac{bP}{p}(1+aq)w \Delta Td\mho-\int_{\mho}Q_{1}\Delta T d\mho\nonumber\\
&=&l_{1}+l_{2}+l_{3}+l_{4}.
\end{eqnarray}
To estimates $l_{1},$ using H$\mathrm{\ddot{o}}$lder inequality, interpolation inequality and Young's inequality we have
\begin{eqnarray*}
l_{1}&\leq& | \Delta T|_{2}|\nabla T|_{4}|v|_{4}\nonumber\\
&\leq & c| \Delta T|_{2}[|\nabla T|_{2}^{\frac{1}{4}}(|\nabla T|_{2}^{\frac{3}{4}}+ |\Delta T|_{2}^{\frac{3}{4}}+|\nabla T_{\xi} |_{2}^{\frac{3}{4}}    )    ]|v|_{4}\nonumber\\
&\leq & \varepsilon |\Delta T|_{2}^{2}+  \varepsilon|\nabla T_{\xi} |_{2}^{2} +c |\nabla T|_{2}^{2}|v|_{4}^{2}+ c|\nabla T|_{2}^{2}|v|_{4}^{8}\\
&\leq & \varepsilon |\Delta T|_{2}^{2}+  \varepsilon|\nabla T_{\xi} |_{2}^{2} +c |\nabla T|_{2}^{2}.
\end{eqnarray*}
To estimates $l_{2},$ using H$\mathrm{\ddot{o}}$lder inequality, Minkowski inequality, interpolation inequality and Young's inequality we obtain
\begin{eqnarray*}
l_{2}&\leq&\int_{S^{2}}\Big{(}\int_{0}^{1}|\mathrm{div} v|d\xi \int_{0}^{1}|T_{\xi}||\Delta T |d\xi   \Big{)}dS^{2}\nonumber\\
&\leq& \int_{S^{2}}\Big{(}\int_{0}^{1}|\mathrm{div} v|d\xi ( \int_{0}^{1}T_{\xi}^{2}d\xi )^{\frac{1}{2}}(\int_{0}^{1} |\Delta T |^{2}d\xi)^{\frac{1}{2}}\Big{)}dS^{2}\nonumber\\
&\leq& |\Delta T|_{2}\Big{(}\int_{S^{2}}(\int_{0}^{1}T_{\xi}^{2}d\xi )^{2}  dS^{2}  \Big{)}^{\frac{1}{4}}\Big{(}\int_{S^{2}}(\int_{0}^{1}|\mathrm{div}  v|d\xi )^{4}dS^{2}\Big{)}^{\frac{1}{4}}\nonumber\\
&\leq& |\Delta T|_{2}\Big{(}\int_{0}^{1}(\int_{S^{2}}T_{\xi}^{4}dS^{2} )^{\frac{1}{2}}d\xi    \Big{)}^{\frac{1}{2}}
\Big{(}\int_{0}^{1}(\int_{S^{2}}|\mathrm{div} v|^{4}dS^{2})^{\frac{1}{4}}d\xi\Big{)}\nonumber\\
&\leq&c |\Delta T|_{2}\Big{(} \int_{0}^{1}|T_{\xi}|_{2}(|T_{\xi}|_{2}+ |\nabla T_{\xi}|_{2}   )    d\xi\Big{)}^{\frac{1}{2}}\nonumber\\
&&\cdot \Big{(}  \int_{0}^{1}|\mathrm{div} v    |_{2}^{\frac{1}{2}}(|\mathrm{div} v    |_{2}^{\frac{1}{2}}+|\Delta v|_{2}^{\frac{1}{2}}   ) d\xi    \Big{)}\nonumber\\
&\leq&c|\Delta T|_{2}|T_{\xi}|_{2}^{\frac{1}{2}}(|T_{\xi}|_{2}^{\frac{1}{2}}+   |\nabla T_{\xi}|_{2}^{\frac{1}{2}}  )\nonumber\\
&&\cdot |\mathrm{div} v|_{2}^{\frac{1}{2}}( |\mathrm{div} v|_{2}^{\frac{1}{2}}+ |\Delta v|_{2}^{\frac{1}{2}})\nonumber\\
&\leq&\varepsilon |\Delta T|_{2}^{2}+c(|T_{\xi}|_{2}^{2}+|T_{\xi}|_{2}|\nabla T_{\xi}|_{2})(|\mathrm{div} v|_{2}^{2}+|\mathrm{div} v|_{2}|\Delta v|_{2} )\\
&\leq&\varepsilon |\Delta T|_{2}^{2}+c(1+|\nabla T_{\xi}|_{2})(1+|\Delta v|_{2} ).
\end{eqnarray*}
Taking an analogous argument as above we get
\begin{eqnarray*}
l_{3}&\leq& \varepsilon |\Delta T|_{2}^{2}+c(1+|\nabla q|_{2})(1+|\Delta v|_{2} ).
\end{eqnarray*}
Similarly, we have
\begin{eqnarray*}
l_{4}\leq \varepsilon |\Delta T |_{2}^{2}+c|Q_{1}|_{2}^{2}.
\end{eqnarray*}
Combining $(3.70)$ and estimates of $l_{1}-l_{4}$ we get
\begin{eqnarray}
&&\frac{1}{2}\frac{d |\nabla T|_{2}^{2}}{dt}+|\Delta T|_{2}^{2}+|\nabla T_{\xi}|_{2}^{2}+\alpha_{s}|\nabla T(\xi=1) |_{2}^{2}\nonumber\\
&\leq& c|\nabla T|_{2}^{2}
+c|Q_{1}|_{2}^{2}+c(1+|\nabla T_{\xi}|_{2})(1  +  |\Delta v |_{2} )\nonumber\\
&&+c(1+|\nabla q|_{2})(1 +   |\Delta v |_{2} ).
\end{eqnarray}
According to $(3.27), (3.33), (3.65 ), (3.69)$ and uniform Gronwall inequality, we have for all $t\geq 0$
\begin{eqnarray}
|\nabla T (t)|_{2}^{2}+ \int_{t}^{t+1}( |\Delta T(t)|_{2}^{2}+|\nabla T_{\xi}(t)|_{2}^{2})dt \leq c,
\end{eqnarray}
where $c$ is a positive constant independent of $t.$
Analogously to the deduction of $(3.71)$, we have
\begin{eqnarray}
&&\frac{1}{2}\frac{d |\nabla q|_{2}^{2}  }{dt}+|\Delta q|_{2}^{2}+|\nabla q_{\xi}|_{2}^{2}+\beta_{s}| \nabla q(\xi=1) |_{2}^{2}\nonumber\\
&\leq& \varepsilon|\Delta q|_{2}^{2} +c|Q_{2}|_{2}^{2}+c|\nabla q|_{2}^{2}\nonumber\\
&&+c(1+|\nabla T_{\xi}|_{2} )(1+|\Delta v|_{2}  ).
\end{eqnarray}
According to $(3.27), (3.65 ), (3.72 )$ and uniform Gronwall inequality, we have for all $t\geq 0$
\begin{eqnarray}
|\nabla q (t)|_{2}^{2}\leq c,
\end{eqnarray}
where $c$ is a positive constant independent of $t.$
\par
Combining the above results, we have proved the existence of an absorbing ball for the strong solution $(v, T, q)$ in the solution space $V$ for the 3D viscous PEs of large-scale moist atmosphere.

\section{Continuity of strong solution with respect to $t$.}
In this section, we will consider the continuity of strong solution with respect to time $t$ in $V$, which will be helpful for proving the existence and connectedness of the global attractor for moist primitive equation. To establish our result of this section, we need the following lemma which is a consequence of a general result of $\cite{LM}.$ For the proof, one can refer to $\cite{T}.$
\begin{lemma}
Let $V,H , V'$ be three Hilbert spaces such that $V\subset H=H'\subset V',$ where $H'$ and $V'$ are the dual spaces of $H$ and $V$ respectively. Suppose $u\in L^{2}([0,T]; V)$ and $u'\in L^{2}([0,T]; V') .$ Then $u$ is almost everywhere equal to a function continuous from $[0,T]$ into $H.$
\end{lemma}
From $(1.11)$ and Lemma $2.1,$ we obtain for $\eta\in V_{1}=D(A_{1}^{\frac{1}{2}}),$
\begin{eqnarray*}
\langle \partial_{t} A_{1}^{\frac{1}{2}}v,      \eta \rangle&=& \langle \partial_{t} v,    A_{1}^{\frac{1}{2}}\eta \rangle \nonumber\\
&=&-\langle  A_{1}v,    A_{1}^{\frac{1}{2}} \eta\rangle-\langle \nabla_{v}v,     A_{1}^{\frac{1}{2}} \eta \rangle\nonumber\\
&&-\langle \Big{(}\int_{\xi}^{1}\mathrm{div} v(x,y,\xi,t )d\xi   \Big{)}v_{\xi},   A_{1}^{\frac{1}{2}} \eta \rangle\nonumber\\
&&-\langle  \int_{\xi}^{1}\frac{bP}{p}\mathrm{grad }[(1+aq)T]d\xi',   A_{1}^{\frac{1}{2}} \eta  \rangle\nonumber\\
&&-\langle \frac{f}{R_{0}} v^{\bot},       A_{1}^{\frac{1}{2}} \eta \rangle.
\end{eqnarray*}
Using H$\mathrm{\ddot{o}}$lder inequality, Agmon inequality and interpolation inequality yields
\begin{eqnarray*}
|\langle A_{1}v,       A_{1}^{\frac{1}{2}}\eta  \rangle|&\leq& |A_{1}v|_{2}\|\eta\|_{1},\\
| \langle (v\cdot \nabla  )  v,      A_{1}^{\frac{1}{2}}\eta  \rangle |&\leq& c |v|_{\infty}\|v\|_{1}| A_{1}^{\frac{1}{2}}\eta |_{2}\leq c\|v\|_{1}^{\frac{3}{2}}\|v\|_{2}^{\frac{1}{2}}| A_{1}^{\frac{1}{2}}\eta |_{2}.
\end{eqnarray*}
By H$\mathrm{\ddot{o}}$lder inequality, Minkowski inequality and interpolation inequality, we have
\begin{eqnarray*}
&&-\langle  \int_{\xi}^{1}\frac{bP}{p}\mathrm{grad }[(1+aq)T]d\xi',   A_{1}^{\frac{1}{2}} \eta  \rangle\\
&\leq&c\int_{S^{2}}\Big{(}\int_{0}^{1}(|\mathrm{grad} T|+ | \mathrm{grad} (Tq) |  )d\xi\int_{0}^{1}A_{1}^{\frac{1}{2}}\eta d\xi \Big{)}dS^{2}\\
&\leq&c \|T\|_{1}|A_{1}^{\frac{1}{2}}\eta|_{2}+c|A_{1}^{\frac{1}{2}}\eta|_{2} \int_{0}^{1}|\mathrm{grad} T|_{4}|q|_{4}d\xi\\
&&+c|A_{1}^{\frac{1}{2}}\eta|_{2} \int_{0}^{1}|\mathrm{grad} q|_{4}|T|_{4}d\xi\\
&\leq&c \|T\|_{1}|A_{1}^{\frac{1}{2}}\eta|_{2}+c|A_{1}^{\frac{1}{2}}\eta|_{2} \int_{0}^{1}|\mathrm{grad} T|_{2}^{\frac{1}{2}}\\
&&\cdot(|\Delta T|_{2}^{\frac{1}{2}}+|\mathrm{grad} T|_{2}^{\frac{1}{2}})
|q|_{4}d\xi\\
&&+c|A_{1}^{\frac{1}{2}}\eta|_{2} \int_{0}^{1}|\mathrm{grad} q|_{2}^{\frac{1}{2}}\\
&&\cdot(|\Delta q|_{2}^{\frac{1}{2}}+|\mathrm{grad} q|_{2}^{\frac{1}{2}})|T|_{4}d\xi\\
&\leq&c \|T\|_{1}|A_{1}^{\frac{1}{2}}\eta|_{2}+c|A_{1}^{\frac{1}{2}}\eta|_{2}\|T\|_{1}^{\frac{1}{2}}\|T\|_{2}^{\frac{1}{2}}|q|_{4}\\
&&+c|A_{1}^{\frac{1}{2}}\eta|_{2}\|q\|_{1}^{\frac{1}{2}}\|q\|_{2}^{\frac{1}{2}}|T|_{4}.
\end{eqnarray*}
Similarly, we have
\begin{eqnarray*}
|\langle \Big{(}\int_{\xi}^{1}\mathrm{div} v(x,y,\xi',t )d\xi'   \Big{)}v_{\xi},   A_{1}^{\frac{1}{2}} \eta \rangle|\leq c\|v\|_{1}\|v\|_{2}\|\eta\|_{1}.
\end{eqnarray*}
In view of the above arguments, we conclude that
\begin{eqnarray*}
\|\partial_{t}(A_{1}^{\frac{1}{2}}v )\|_{V_{1}'}&\leq& c\|v\|_{2}+c\|v\|_{1}\|v\|_{2}+c\|T\|_{1}\\
&&+c(\|T\|_{1}^{\frac{1}{2}} \|T\|_{2}^{\frac{1}{2}}|q|_{4}+ \|q\|_{1}^{\frac{1}{2}} \|q\|_{2}^{\frac{1}{2}}|T|_{4}   ).
\end{eqnarray*}
Since
\begin{eqnarray*}
v \in L^{\infty} ([0, \tau]; V_{1})\cap L^{2}([0,\tau]; H^{2}(\mho)),\ \ \ \forall \ \tau>0,
\end{eqnarray*}
we have
\begin{eqnarray*}
A_{1}^{\frac{1}{2}}v\in L^{2}([0,\tau]; V_{1}),\ \ \ \ \   \partial_{t}A_{1}^{\frac{1}{2}}v \in L^{2}([0,\tau]; V_{1}'), \ \ \ \forall \ \tau>0.
\end{eqnarray*}
Therefore, by Lemma $4.1$, we infer that $v\in C([0,\tau]; V_{1}).$
For any $\phi \in V_{2}=D(A_{2}^{\frac{1}{2}}),$ we have
\begin{eqnarray*}
\langle \partial_{t}A_{2}^{\frac{1}{2}}T,      \phi \rangle&=& \langle\partial_{t} T,    A_{2}^{\frac{1}{2}} \phi \rangle\\
&=&-\langle A_{2}T,   A_{2}^{\frac{1}{2}} \phi \rangle-\langle \nabla_{v}T,   A_{2}^{\frac{1}{2}} \phi  \rangle\\
&&-\langle \Big{(}\int_{\xi}^{1}\mathrm{div} v(x,y,\xi',t )d\xi'   \Big{)}\partial_{\xi}T,    A_{2}^{\frac{1}{2}} \phi \rangle\\
&&+\langle \frac{bP}{p}(1+aq) \Big{(}\int_{\xi}^{1}\mathrm{div} v(x,y,\xi',t )d\xi'   \Big{)},     A_{2}^{\frac{1}{2}} \phi \rangle+\langle Q_{1},   A_{2}^{\frac{1}{2}} \phi \rangle.
\end{eqnarray*}
By H$\mathrm{\ddot{o}}$lder inequality, Minkowski inequality and interpolation inequality, we have
\begin{eqnarray*}
&&\langle \frac{bP}{p}(1+aq)\Big{(}\int_{\xi}^{1}\mathrm{div} v(x,y,\xi',t )d\xi'   \Big{)},  A_{2}^{\frac{1}{2}}\phi  \rangle\\
&\leq& c\|v\|_{1}\|\phi\|_{1}+c\int_{S^{2}}\Big{(}\int_{0}^{1}|\mathrm{div} v|d\xi\Big{)} \Big{(}\int_{0}^{1}|q| |A_{1}^{\frac{1}{2}}\phi|d\xi \Big{)}dS^{2}\\
&\leq&c\|v\|_{1}\|\phi\|_{1}+c\int_{S^{2}}\Big{(}\int_{0}^{1}|\mathrm{div} v|d\xi\Big{)}\Big{(}\int_{0}^{1}q^{2}d\xi\Big{)}^{\frac{1}{2}}
\Big{(}\int_{0}^{1}|A_{2}^{\frac{1}{2}}\phi|^{2}d\xi\Big{)}^{\frac{1}{2}}dS^{2}\\
&\leq&c\|v\|_{1}\|\phi\|_{1}+c\|\phi\|_{1} \Big{(}\int_{0}^{1}|\mathrm{div} v|_{4}d\xi\Big{)}  \Big{(}\int_{0}^{1}|q|_{4}^{2}d\xi\Big{)}^{\frac{1}{2}}\\
&\leq&c\|v\|_{1}\|\phi\|_{1}+c\|\phi\|_{1}|q|_{4}\|v\|_{1}^{\frac{1}{2}}\|v\|_{2}^{\frac{1}{2}}.
\end{eqnarray*}
Consequently, taking an analogous argument about $\|\partial_{t} A_{1}^{\frac{1}{2}}v\|_{V_{1}'}$, we have
\begin{eqnarray*}
\|\partial_{t}A_{2}^{\frac{1}{2}}T\|_{V_{2}'}&\leq& c(\|T\|_{2}+\|v\|_{1}+\|v\|_{1}^{\frac{1}{2}}\|v\|_{2}^{\frac{1}{2}}\|T\|_{1}+ |q|_{4}\|v\|_{1}^{\frac{1}{2}}\|v\|_{2}^{\frac{1}{2}}
\\
&&+\|T\|_{1}^{\frac{1}{2}}\|T\|_{2}^{\frac{1}{2}}\|v\|_{1}^{\frac{1}{2}}\|v\|_{2}^{\frac{1}{2}}+|Q_{1}|_{2}    ).
\end{eqnarray*}
Similarly, we can infer that
\begin{eqnarray*}
\|\partial_{t}A_{3}^{\frac{1}{2}}q\|_{V_{2}'}\leq c(\|q\|_{2}+\|v\|_{1}^{\frac{1}{2}}\|v\|_{2}^{\frac{1}{2}}\|q\|_{1}+\|q\|_{1}^{\frac{1}{2}}\|q\|_{2}^{\frac{1}{2}}\|v\|_{1}^{\frac{1}{2}}\|v\|_{2}^{\frac{1}{2}}
+|Q_{2}|_{2}    ).
\end{eqnarray*}
Since
\begin{eqnarray*}
(v, T, q) \in L^{\infty} ([0, \tau]; V)\cap L^{2}([0,\tau]; H^{2}(\mho)),\ \ \ \forall \ \tau>0,
\end{eqnarray*}
we have
\begin{eqnarray*}
A_{2}^{\frac{1}{2}}T\in L^{2}([0,\tau]; V_{2}),\ \ \ \ \   \partial_{t}A_{2}^{\frac{1}{2}}T \in L^{2}([0,\tau]; V_{2}')
\end{eqnarray*}
and
\begin{eqnarray*}
A_{3}^{\frac{1}{2}}q\in L^{2}([0,\tau]; V_{3}),\ \ \ \ \   \partial_{t}A_{3}^{\frac{1}{2}}q \in L^{2}([0,\tau]; V_{3}').
\end{eqnarray*}
By Lemma $4.1,$ we infer that $T,q\in C([0, +\infty); V_{2})$ and $ C([0, +\infty); V_{3}),$ respectively. So far, we obtain that
\begin{eqnarray*}
(v,T,q)\in C([0, +\infty); V).
\end{eqnarray*}
\section{Continuity in $V$ with respect to initial data.}
It is shown in $\cite{GH2}$ that the  strong solution to the 3D viscous primitive equations of large-scale moist atmosphere is unique and Lipschitz continuous with respect to the initial condition in $H.$ But what we need to do here is to show the continuity property in $V.$

In the following, will prove that for any fixed $t>0,$ the mapping $(v_{0}, T_{0}, q_{0})\mapsto (v(t), T(t), q(t) )$ is Lipschitz continuous from $V$ into itself for all the strong solutions.

Assume $(v_{i}, T_{i}, q_{i} ), i=1,2,$ are two strong solutions to the equation $(1.11)-(1.17)$ with initial data $(v_{0,i}, T_{0,i}, q_{0,i}  )\in V.$ Let
\begin{eqnarray*}
u=v_{1}-v_{2},\ \ \tau=T_{1}-T_{2},\ \ q=q_{1}-q_{2},\ \ \Phi_{s}(t; \theta, \varphi )= \Phi_{1,s}(t; \theta, \varphi )-\Phi_{2,s}(t; \theta, \varphi ).
\end{eqnarray*}
Then we derive from $(1.11)-(1.17)$ that
\begin{eqnarray}
\partial_{t}u&+&L_{1}u+\nabla_{v_{1}}u+\nabla_{u}v_{2}+\Big{(}\int_{\xi}^{1}\mathrm{div} v_{1}(x,y,\xi',t )d\xi'   \Big{)}\partial_{\xi}v\nonumber\\
&&+\Big{(}\int_{\xi}^{1}\mathrm{div} v(x,y,\xi',t )d\xi'   \Big{)}\partial_{\xi}v_{2}+\frac{f}{R_{0}}v^{\bot}+ \mathrm{grad} \Phi_{s}\nonumber\\
&&+\int_{\xi}^{1}\frac{bP}{p} \mathrm{grad} T d\xi'+\int_{\xi}^{1}\frac{abP}{p} \mathrm{grad} (q_{1} T) d\xi'+\int_{\xi}^{1}\frac{abP}{p} \mathrm{grad} (q T_{2}) d\xi'=0,\\
\partial_{t}T&+&L_{2}T+\nabla_{v_{1}}T+\nabla_{v}T_{2}+\Big{(}\int_{\xi}^{1}\mathrm{div} v_{1}(x,y,\xi',t )d\xi'   \Big{)}\partial_{\xi}T\nonumber\\
&&+\Big{(}\int_{\xi}^{1}\mathrm{div} v(x,y,\xi',t )d\xi'   \Big{)}\partial_{\xi}T_{2}-\frac{bP}{p}\Big{(}\int_{\xi}^{1}\mathrm{div} v(x,y,\xi',t )d\xi'   \Big{)}\nonumber\\
&&-\frac{abP}{p}q_{1}\Big{(}\int_{\xi}^{1}\mathrm{div} v(x,y,\xi',t )d\xi'   \Big{)}
-\frac{abP}{p}q\Big{(}\int_{\xi}^{1}\mathrm{div} v_{2}(x,y,\xi',t )d\xi'   \Big{)}=0,\\
\partial_{t}q&+&L_{3}q+\nabla_{v_{1}}q+\nabla_{v}q_{2}+\Big{(}\int_{\xi}^{1}\mathrm{div} v_{1}(x,y,\xi',t )d\xi'   \Big{)}\partial_{\xi}q \nonumber\\ &&+\Big{(}\int_{\xi}^{1}\mathrm{div} v(x,y,\xi',t )d\xi'   \Big{)}\partial_{\xi}q_{2}=0,
\end{eqnarray}
\begin{eqnarray}
u|_{t=0}= v_{0}^{1}-v_{0}^{2},\ \  T|_{t=0}=T_{0}^{1}-T_{0}^{2},\ \ q|_{t=0}=q_{0}^{1}-q_{0}^{2},
\end{eqnarray}
\begin{eqnarray}
\xi=1:\ \ \partial_{\xi}v=0,\ \ \partial_{\xi}T=-\alpha_{s}T,\ \ \partial_{\xi} q=-\beta_{s}q,
\end{eqnarray}
\begin{eqnarray}
\xi=0:\ \ \partial_{\xi}v=0,\ \ \partial_{\xi}T=0,\ \ \partial_{\xi} q=0.
\end{eqnarray}
Taking inner product of $(5.75)$ with $A_{1}u$ in $L^{2}(\mho)$ we obtain
\begin{eqnarray}
\frac{1}{2}\frac{d}{dt}\|u\|_{1}^{2}+|A_{1}u|_{2}^{2}&=& -\langle \nabla_{v_{1}}u,     A_{1}u \rangle-\langle \nabla_{u}v_{2},    A_{1}u \rangle\nonumber\\
&&-\langle\Big{(}\int_{\xi}^{1}\mathrm{div} u(x,y,\xi',t )d\xi'   \Big{)}\partial_{\xi}v_{2},    A_{1} u  \rangle\nonumber\\
&&-\langle \Big{(}\int_{\xi}^{1}\mathrm{div} v_{1}(x,y,\xi',t )d\xi'   \Big{)}\partial_{\xi} u,   A_{1}u \rangle\nonumber\\
&&-\langle \int_{\xi}^{1}\frac{bP}{p} \mathrm{grad }T d\xi',    A_{1} u  \rangle-\langle (\frac{f}{R_{0}} u^{\bot}+\mathrm{grad} \Phi_{s}),   A_{1}u   \rangle\nonumber\\
&&-\langle \int_{\xi}^{1}\frac{abP}{p} \mathrm{grad}(q_{1}T)d\xi',     A_{1} u    \rangle-\langle \int_{\xi}^{1}\frac{abP}{p} \mathrm{grad}(qT_{2})d\xi'   ,     A_{1} u\rangle\nonumber\\
&=&\sum_{i=1}^{8}k_{i}.
\end{eqnarray}
By  H$\mathrm{\ddot{o}}$lder inequality, Agmon inequality and Young's inequality , we have
\begin{eqnarray*}
k_{1}&\leq& |v_{1}|_{\infty}(|\nabla_{e_{\theta}}u|_{2}+  |\nabla_{e_{\varphi}}u|_{2})|A_{1}u|_{2}\\
&\leq &c\|v_{1}\|_{1}^{\frac{1}{2}}|A_{1}v_{1}|_{2}^{\frac{1}{2}}\|u\|_{1}|A_{1}u|_{2}\\
&\leq & \varepsilon |A_{1}u|_{2}^{2}+c\|v_{1}\|_{1}|A_{1}v_{1}|_{2}\|u\|_{1}^{2}.
\end{eqnarray*}
Similarly, we obtain
\begin{eqnarray*}
k_{2}\leq |u|_{\infty}\|v_{2}\|_{1}|A_{1}u|_{2}\leq \varepsilon |A_{1}u|_{2}^{2}+c\|u\|_{1}^{2}\|v_{2}\|_{1}^{4}.
\end{eqnarray*}
Taking an analogous argument as $(3.61),$ we have
\begin{eqnarray*}
k_{3}+k_{4}\leq \varepsilon |A_{1}u|_{2}^{2}+c\|u\|_{1}^{2}\|v_{2}\|_{1}^{2}\| v_{2}\|_{2}^{2}+c\|u\|_{1}^{2}\|v_{1}\|_{1}^{2}\|v_{1}\|_{2}^{2} .
\end{eqnarray*}
In view of Lemma $2.1$ and H$\mathrm{\ddot{o}}$ler inequality, we obtain
\begin{eqnarray*}
k_{5}+k_{6}\leq \varepsilon |A_{1}u|_{2}^{2}+c\|T\|_{1}^{2}+c|u|_{2}^{2}.
\end{eqnarray*}
To estimate $k_{7}$,  by H$\mathrm{\ddot{o}}$lder inequality, Minkowski inequality and interpolation inequality  we have
\begin{eqnarray*}
k_{7}&\leq& c\int_{S^{2}}(\int_{0}^{1}|\mathrm{grad} q_{1}| |T|d\xi \int_{0}^{1}|A_{1}u|d\xi  )dS^{2}\\
&&+c\int_{S^{2}}(\int_{0}^{1}| q_{1}| |\mathrm{grad} T|d\xi \int_{0}^{1}|A_{1}u|d\xi  )dS^{2}\\
&\leq &c |A_{1}u|_{2} \Big{(}\int_{S^{2}}(\int_{0}^{1}|\mathrm{grad} q_{1}|^{2}d\xi)^{2}dS^{2}\Big{)}^{\frac{1}{4}}\Big{(}\int_{S^{2}}(\int_{0}^{1}|T|^{2}d\xi)^{2}dS^{2}\Big{)}^{\frac{1}{4}}\\
&&+c |A_{1}u|_{2} \Big{(}\int_{S^{2}}(\int_{0}^{1}| q_{1}|^{2}d\xi)^{2}dS^{2}\Big{)}^{\frac{1}{4}}\Big{(}\int_{S^{2}}(\int_{0}^{1}|\mathrm{grad} T|^{2}d\xi)^{2}dS^{2}\Big{)}^{\frac{1}{4}}\\
&\leq & \varepsilon |A_{1}u|_{2}^{2}+c|\mathrm{grad }q_{1}|_{2}( |\mathrm{grad }q_{1}|_{2}+|\Delta q_{1}|_{2}   )(|T|_{2}^{2}+|T|_{2}|\nabla T|_{2})\\
&&+c|\mathrm{grad }T|_{2}( |\mathrm{grad }T|_{2}+|\Delta T|_{2}   )|q_{1}|_{4}^{2}\\
&\leq & \varepsilon |A_{1}u|_{2}^{2}+\varepsilon |\Delta T|_{2}^{2}+c\|T\|_{1}^{2}(\|q_{1}\|_{1}^{4}+\|q_{1}\|_{2}^{2} ).
\end{eqnarray*}
Similarly, we have
\begin{eqnarray*}
k_{8}
&\leq& \varepsilon |A_{1}u|_{2}^{2}  +c|\mathrm{grad }q|_{2}( |\mathrm{grad }q|_{2}+|\Delta q|_{2}   )|T_{2}|_{4}^{2}\\
&&+c|\mathrm{grad }T_{2}|_{2}( |\mathrm{grad }T_{2}|_{2}+|\Delta T_{2}|_{2}   )(|q|_{2}^{2}+|q|_{2}|\nabla q|_{2})\\
&\leq & \varepsilon |A_{1}u|_{2}^{2}+\varepsilon |\Delta q|_{2}^{2}+c\|q\|_{1}^{2}(\|T_{2}\|_{1}^{4}+\|T_{2}\|_{2}^{2} ).
\end{eqnarray*}
By $(5.81)$ and estimates of $k_{1}-k_{8},$  we get
\begin{eqnarray}
\frac{1}{2}\frac{d\|u\|_{1}^{2}}{dt}+|A_{1}u|_{2}^{2}&\leq& \varepsilon|A_{1}u|_{2}^{2}+\varepsilon|A_{2}T|_{2}^{2}+\varepsilon|A_{3}q|_{2}^{2}\nonumber\\
&&+c\|u\|_{1}^{2}(\|v_{1}\|_{1}^{2}\|v_{1}\|_{2}^{2}+ \|v_{2}\|_{1}^{2}\|v_{2}\|_{2}^{2}+\|v_{2}\|_{1}^{4}+1)\nonumber\\
&&+c\|T\|_{1}^{2}(\|q_{1}\|_{1}^{4}+\|q_{1}\|_{2}^{2})\nonumber\\
&&+c\|q\|_{1}^{2}(\|T_{2}\|_{1}^{4}+\|T_{2}\|_{2}^{2}).
\end{eqnarray}
Taking an analogous argument as above, from $(5.76)$ and $(5.77)$ we have
\begin{eqnarray}
&&\frac{1}{2}\frac{d\| T\|_{1}^{2}}{dt}+|A_{2}T|_{2}^{2}\nonumber\\
&\leq&\varepsilon |A_{2}T|_{2}^{2}+\varepsilon |A_{1}u|_{2}^{2}+c\|q\|_{1}^{2}\|v_{2}\|_{1}\|v_{2}\|_{2}\nonumber\\
&&+c\|T\|_{1}^{2}(1+\|v_{1}\|_{1}^{2}\|v_{1}\|_{2}^{2} )\nonumber\\
&&+c \|u\|_{1}^{2}(1+\|q_{1}\|_{1}^{4}+\|T_{2}\|_{2}^{2}+\|T_{2}\|_{1}^{2}\|T_{2}\|_{2}^{2}  )
\end{eqnarray}
and
\begin{eqnarray}
&&\frac{1}{2}\frac{d\|q\|_{1}^{2}}{dt}+|A_{3}q|_{2}^{2}\nonumber\\
&\leq&\varepsilon |A_{3}q|_{2}^{2}+c\|q\|_{1}^{2}(1+ \|v_{1}\|_{1}^{2}\|v_{1}\|_{2}^{2})\nonumber\\
&&+c\|u\|_{1}^{2}(  \|q_{2}\|_{2}^{2}+  \|q_{2}\|_{1}^{2}\|q_{2}\|_{2}^{2}).
\end{eqnarray}
Let
\begin{eqnarray*}
g_{1}&:=&1+\|v_{1}\|_{1}^{2}\|v_{1}\|_{2}^{2}+\|v_{2}\|_{1}^{4}+\|v_{2}\|_{1}^{2}\|v_{2}\|_{2}^{2}+\|T_{2}\|_{2}^{2}\nonumber\\
&&+\|T_{2}\|_{1}^{2}\|T_{2}\|_{2}^{2}+\|q_{1}\|_{1}^{4}+\|q_{2}\|_{2}^{2}+\|q_{2}\|_{1}^{2}\|q_{2}\|_{2}^{2},
\end{eqnarray*}
\begin{eqnarray*}
g_{2}:=1+\|q_{1}\|_{2}^{2}+\|q_{1}\|_{1}^{4}+\|v_{1}\|_{1}^{2}\|v_{1}\|_{2}^{2}
\end{eqnarray*}
and
\begin{eqnarray*}
g_{3}:=1+\|T_{2}\|_{2}^{2}+\|T_{2}\|_{1}^{4}+\|v_{2}\|_{1}\|v_{2}\|_{2} +\|v_{1}\|_{1}^{2}\|v_{1}\|_{2}^{2}.
\end{eqnarray*}
Obviously, for arbitrary $0\leq a< b<\infty$, we have
\begin{eqnarray*}
\int_{a}^{b}(g_{1}(t)+g_{2}(t)+g_{3}(t))dt< \infty.
\end{eqnarray*}
Thereby, in view of $(5.82)-(5.84)$ we get
\begin{eqnarray*}
\frac{d(\|u\|_{1}^{2}+\|T\|_{1}^{2}+\|q\|_{1}^{2}) }{dt}\leq c( g_{1}(t)+g_{2}(t)+g_{3}(t) )(\|u\|_{1}^{2}+\|T\|_{1}^{2}+\|q\|_{1}^{2}),
\end{eqnarray*}
 which combined with Gronwall lemma implies
\begin{eqnarray*}
&&\|u(t)\|_{1}^{2}+\|T(t)\|_{1}^{2}+\|q(t)\|_{1}^{2}\\
&\leq& c(\|v_{0,1}-v_{0,2}\|_{1}^{2} + \|T_{0,1}-T_{0,2}\|_{1}^{2} +\|q_{0,1}-q_{0,2}\|_{1}^{2}  )
e^{\int_{0}^{t}(g_{1}(s)+g_{2}(s)+g_{3}(s))ds }.
\end{eqnarray*}
So far, we have shown that for $t>0,\ (v(t), T(t), q(t) )$ is Lipschitz continuous in $V$ with respect to the initial data $(v(0), T(0), q(0) ) .$

\section{The global attractor.}
In this section, we present our main result, the existence of global attractor, of this paper. To show our main result, we make use of the following theorem from Teman $\cite{T}.$ For more details about the theorem, we can see, e.g.,  $\cite{BV, CV, H, L, R, SY, T}$ and the references therein.
\begin{theorem}
Suppose that $X$ is a metric space and semigroup $\{ S(t)\}_{t\geq 0}$ is a family of operators from $X$ into itself such that

$(i)$ for any fixed $t>0, S(t)$ is continuous from $X$ into itself;

$(ii)$ for some $t_{0}, S(t_{0})$ is compact from $X$ into itself;

$(iii)$ there exists a subset $B_{0}$ of $X$ which is bounded, a subset $U$ of $X$ is open, such that

\qquad $B_{0}\subseteq U\subseteq X,$ and $B_{0}$ is the absorbing set of $U,$ i.e. for any bounded subset $B\subset U,$ there

\qquad is a $t_{0}=t_{0}(B),$ such that
\begin{eqnarray*}
S(t)B\subset B_{0}, \ \ \ \  \ \forall\ t>t_{0}(B).
\end{eqnarray*}
Then $\mathcal{A}:=\omega(B_{0}),$ the $\omega-$limit set of $B_{0}$, is a compact attractor which attracts all the bounded sets of $U,$ i.e. for any $x\in U,$
\begin{eqnarray*}
\lim\limits_{t\rightarrow \infty} dist ( S(t)x,  \mathcal{A  })=0.
\end{eqnarray*}
The set $\mathcal{A} $ is the maximal bounded attractor in $U$ for the inclusion relation.\\

Suppose in addition that $X$ is a Banach space, $U$ is a convex and\\

$(iv)\ \forall x\in X, S(t)x: \mathbb{R}_{+}\mapsto X$ is continuous.\\
Then $\mathcal{A}:=\omega(B_{0})  $ is also connected.\\

If $U=X, \mathcal{A}$ is the global attractor of the semigroup $\{S(t)   \}_{t\geq 0}$ in $X.$
\end{theorem}
Next, we give our main result of the paper and complete the proof of the theorem.
\begin{theorem}
Assume $Q_{1}, Q_{2}\in L^{2}(\mho)$ and $ Q_{1}|_{\xi=1}, Q_{2}|_{\xi=1}\in L^{2}(S^{2}).$ Then, for $t\geq 0,$ the solution operator $\{S(t)   \}_{t\geq 0}$ of the $3D$ viscous PEs of large-scale moist atmosphere $(1.11)-(1.17): S(t)(v_{0}, T_{0}, q_{0})= (v(t), T(t), q(t) )$ defines a semigroup in the space $V.$ Furthermore, the results below hold:
\par
$(1)$ For any $(v_{0}, T_{0}, q_{0})\in V, t\rightarrow S(t)(v_{0}, T_{0}, q_{0}) $ is a continuous map from $\mathbb{R}_{+} $ into $V.$
\par
$(2)$ For any $t>0, S(t)$ is a continuous map in $V.$
\par
$(3)$ For any $t>0, S(t)$ is a compact map in $V.$
\par
$(4)$ $\{S(t)\}_{t\geq 0}$ possesses a global attractor $\mathcal{A}$ in $V.$ The global attractor in $\mathcal{A}$ is compact and

\quad connected in $V$ and is the maximal bounded attractor in $V$ in the sense of set inclusion

\quad  relation; $\mathcal{A }$ attracts all bounded subset in $V$ in the norm of $V.$
\end{theorem}

To prove Theorem 6.2, we need to check the conditions $(i)-(iv)$ in Theorem 6.1. First,  condition $(i)$, the continuous dependence  on initial data of the solution,  is verified in section $5.$  Second, condition $(iii)$, the existence of an absorbing ball in $V,$ is proved in section $3.$  Third, condition $(iv)$, the regularity of the solution, is shown in section $4.$ Finally, only the condition $(ii),$ compactness of the solution operator $\{S(t)\}_{t\geq 0},$ is left to be checked. We will use Aubin-Lions lemma stated below and continuity argument to verify condition $(ii).$ For more details of the lemma, we can see $\cite{A}, \cite{Li}$ and references therein.
\begin{lemma}
Let $\mathcal{H}_{0}, \mathcal{H}, \mathcal{H}_{1}  $ be Banach spaces such that $\mathcal{H}_{0}, \mathcal{H}_{1}$ are reflexive and $\mathcal{H}_{0}\overset{c}\subset \mathcal{H}\subset \mathcal{H}_{1}.$ Define, for $0<\tau< \infty,$  \begin{eqnarray*}
X:=\Big{\{} u\Big{|}  u\in L^{2}([0,\tau]; \mathcal{H}_{0}),\ \ \frac{d u}{dt}\in L^{2}([0,\tau]; \mathcal{H}_{1})\Big{\}}.
\end{eqnarray*}
Then $X$ is a Banach space equipped with the norm $\|u\|_{L^{2}([0,\tau]; \mathcal{H}_{0}) }+ \|u'\|_{L^{2}([0,\tau]; \mathcal{H}_{1}) } .$ Moreover, $X\overset{c}\subset L^{2}([0,\tau]; \mathcal{H}).$
\end{lemma}
Proof of Theorem $6.1.$ By the argument above, we only have to check condition $(ii),$ the compactness of the solution operator. For any fixed $\tau >0,$ let  $\mathcal{B}$ be a bounded subset of $V$ and $A_{\tau}$ denote the subset of the space $L^{2}([0,T]; H):$
\begin{eqnarray*}
A_{\tau}:=\Big{\{}\Big{(} A_{1}^{\frac{1}{2}}v,   A_{2}^{\frac{1}{2}}T,  A_{3}^{\frac{1}{2}}q   \Big{)}   \Big{|}(v_{0}, T_{0}, q_{0})\in \mathcal{B},   (v(t), T(t), q(t))=S(t)(v_{0}, T_{0}, q(0)), t\in [0, \tau]  \Big{\}}.
\end{eqnarray*}
For $(v_{0}, T_{0}, q(0) ) \in \mathcal{B}$, it has been shown previously that the strong solution $(v, T, q)$ satisfies
\begin{eqnarray*}
( A_{1}^{\frac{1}{2}}v,   A_{2}^{\frac{1}{2}}T,  A_{3}^{\frac{1}{2}}q   )\in L^{2}([0, \tau]; V),\ \ (\partial_{t}A_{1}^{\frac{1}{2}}v, \partial_{t}A_{2}^{\frac{1}{2}}T, \partial_{t} A_{3}^{\frac{1}{2}}q  )\in L^{2}([0, \tau]; V').
\end{eqnarray*}
If we denote
\begin{eqnarray*}
\mathcal{H}_{0}=V, \mathcal{H}=H, \mathcal{H}_{1}=V',
\end{eqnarray*}
from Lemma $6.1$ we infer that  $A_{\tau}$ is compact in $L^{2}([0,\tau]; H).$
\par
To prove solution operator $\{S(t)\}_{t\geq 0}$ is compact in $V,$ for any bounded sequence $\{(v_{0,n}, T_{0,n}, q_{0,n}  )  \}_{0}^{n}\subset \mathcal{B},$ we should show there exists a convergent subsequence of $\{S(t)(v_{0,n}, T_{0,n}, q_{0,n} )   \}_{0}^{\infty}$ in $V.$
\par
Since $A_{\tau}$ is compact in $L^{2}([0,\tau]; H),$ there exists a function $(v_{*}, T_{*}, q_{*} )\in L^{2}([0,\tau]; V)$ such that there is a subsequence of $\{S(\cdot)(v_{0,n}, T_{0,n}, q_{0,n} )   \}_{0}^{\infty},$ still denoted as $\{S(\cdot)(v_{0,n}, T_{0,n}, q_{0,n} )   \}_{0}^{\infty}, $ satisfying
\begin{eqnarray*}
\lim\limits_{n\rightarrow \infty}\int_{0}^{T}\|S(t)(v_{0,n}, T_{0,n}, q_{0,n})-(v_{*}(t), T_{*}(t), q_{*}(t))\|_{1}^{2}dt=0,
\end{eqnarray*}
which implies that there is a subsequence of $\{S(\cdot)(v_{0,n}, T_{0,n}, q_{0,n} )   \}_{0}^{\infty},$ still denoted as $ \{S(\cdot)(v_{0,n}, T_{0,n}, q_{0,n} )   \}_{0}^{\infty}$ for simplicity of notation, converging to $(v_{*}, T_{*}, q_{*}) $ in $V$ a.e. in $(0, \tau):$
\begin{eqnarray*}
\lim\limits_{n\rightarrow \infty}\|S(t)(v_{0,n}, T_{0,n}, q_{0,n})-(v_{*}(t), T_{*}(t), q_{*}(t))\|_{1}=0,\ \ a.e.\ t\ in\ (0, \tau).
\end{eqnarray*}
For any $t\in (0, \tau),$ we can choose a $t_{0}\in (0, t)$ such that
\begin{eqnarray*}
\lim\limits_{n\rightarrow \infty}\|S(t_{0})(v_{0,n}, T_{0,n}, q_{0,n})-(v_{*}(t_{0}), T_{*}(t_{0}), q_{*}(t_{0}))\|_{1}=0.
\end{eqnarray*}
Therefore, for any $t>0,$ by the continuity of $ S(t)$ in $V$, we have
\begin{eqnarray*}
&&\lim\limits_{n\rightarrow \infty}\|S(t)(v_{0,n}, T_{0,n}, q_{0,n})-S(t-t_{0})(v_{*}(t_{0}), T_{*}(t_{0}), q_{*}(t_{0}))\|_{1}\\
&=& \lim\limits_{n\rightarrow \infty}\|S(t-t_{0})S(t_{0})(v_{0,n}, T_{0,n}, q_{0,n})-S(t-t_{0})(v_{*}(t_{0}), T_{*}(t_{0}), q_{*}(t_{0}))\|_{1} =0,
\end{eqnarray*}
which implies $S(t)$ is a compact map in $V$ for any $t>0.$ Then by the above argument of this section and Theorem 6.1, we obtain the existence of the global attractor for the moist primitive equations $(1.11)-(1.17).$
\hspace{\fill}$\square$

\def\refname{ Bibliography}


\begin{thebibliography}{2}

\bibitem {A}
J. Aubin, Un th\'{e}or\`{e}me de compacit\'{e}, C. R. Acad. Sci. Paris, 256 (1963), 5042--5044.

\bibitem {Bje}
V. Bjerknes, Das Problem der Wettervorhersage, betrachtet vom Standpunkte der Mechanik und
der Physik, Meteorol. Z. 21( 1904), 1--7.

\bibitem {BB}
 A.J. Bourgeois, J.T. Beale, Validity of the quasigeostrophic model for large-scale flow in the atmosphere and ocean, SIAM J.
Math. Anal. 25 (1994), 1023--1068.

\bibitem {BV}
A. Babin and M. Vishik, ¡°Attractor of Evolution Equations,¡± North-Holland, Amsterdam
1992.


\bibitem {C}
D. Cordoba, Nonexistence of simple hyperbolic blow-up for the quasi-geostrophic equation, Ann. of Math. 148 (1998),
1135--1152.




\bibitem {CINT}
C. Cao,  S. Ibrahim, K. Nakanishi, E. S. Titi,  Finite-time blowup for the inviscid primitive
equations of oceanic and atmospheric dynamics, Comm. Math. Phys. 337(2015), 473--482.


\bibitem {CLT1}
C. Cao, J. Li, E.S. Titi,  Local and global well-posedness of strong solutions to the 3D primitive equations
with vertical eddy diffusivity, Arch. Anal. Ration. Mech. 214(2014), 35--76.

\bibitem {CLT2}
C. Cao, J. Li, E.S. Titi, Global well-posedness of strong solutions to the 3D primitive equations with
horizontal eddy diffusivity, J. Differ. Equ. 257(2014), 4108--4132.

\bibitem {CLT3}
C. Cao, J. Li, E.S. Titi, Global well-posedness for the 3D primitive equations with only horizontal
viscosity and diffusion, Commun. Pure Appl. Math. Vol.LXIX(2016), 1492--1531.


\bibitem {CMT1}
P. Constantin, A. Majda, E. Tabak, Formation of strong fronts in the 2-D quasigeostrophic thermal active scalar, Nonlinearity 7 (1994), 1495--1533.

\bibitem {CMT2}
P. Constantin, A. Majda, E. Tabak, Singular front formation in a model for quasigeostrophic flow, Phys. Fluids 6 (1994), 9--11.


\bibitem {CT1}
C. Cao  and E. Titi,  Global well-posedness of the three-dimensional viscous primitive equations of
large scale ocean and atmosphere dynamics,  Ann. Math. 166(2007), 245--267.

\bibitem {CT2}
C. Cao, E.S.  Titi,  Global well-posedness of the three-dimensional primitive equations with
partial vertical turbulence mixing heat diffusion, Commun. Math. Phys. 310(2012), 537--568.



\bibitem {CV}
V. Chepyzhov and M. Vishik, ¡°Attractors for Equations of Mathematical Physics,¡± American
Mathematical Society Colloquium Publications, 49. American Mathematical Society, Providence, RI, 2002.

\bibitem {CW}
P. Constantin, J. Wu, Behavior of solutions of 2D quasi-geostrophic equations, SIAM J. Math. Anal. 30 (1999), 937--948.


\bibitem {EM}
P.F. Embid, A.J. Majda, Averaging over fast gravity waves for geophysical flows with arbitrary potential vorticity, Comm.
Partial Differential Equations, 21 (1996), 619--658.


\bibitem {FP}
C. Foias and G. Prodi, Sur le comportement global des solutions non-stationnaires des
\'{e}quations de Navier-Stokes en dimension 2, Rend. Sem. Mat. Univ. Padova 39 (1967), 1--34.


\bibitem {G}
 A. E. Gill,  Atmosphere--Ocean Dynamics (International Geophysics Series vol 30) (San Diego, CA:
Academic), 1982.


\bibitem {GH1}
B. Guo, D. Huang, Existence of weak solutions and trajectory attractors for the moist atmospheric equations in geophysics,
J. Math. Phys. 47(2006), 083508.


\bibitem {GH2}
B. Guo, D. Huang, Existence of the universal attractor for the 3-D viscous primitive equations of large-scale moist atmosphere, J.Differential Equations, 251(2011), 457--491.



\bibitem {GMR}
F. Guill$\acute{e}$n--Gonz$\acute{a}$ez, N. Masmoudi and M.A. Rodr\'{\i}guez--Bellido, Anisotropic estimates and
strong solutions of the Primitive Equations, Diff. Integral Eq. 14(2001), 1381--1408.

\bibitem {H}
J. Hale, ¡°Asymptotic Behavior of Dissipative Systems,¡± American Mathematical Society,
Providence, 1988.

\bibitem {Ha}
G. J. Haltiner,   Numerical Weather Prediction (New York: Wiley), 1971.


\bibitem {HW}
G. J. Haltiner  and R. T. Williams,   Numerical Prediction and Dynamic Meteorology (New York: Wiley), 1980.

\bibitem {Ho}
J.R. Holton, An Introduction to Dynamic Meteorology, third edition, Academic Press, 1992.

\bibitem {HTZ}
C. Hu, R. Temam and M. Ziane, The primitive equations on the large scale ocean under the
small depth hypothesis, Discrete Contin. Dyn. Syst. 9(2003), 97--131.


\bibitem {J}
N.Ju, The global attractor for the solutions to the 3d viscous primitive equations, Discrete and Continuous Dynamical Systems, 17(2007), 159--179.

\bibitem {JT}
N.Ju and R.Teman, Finite dimensions of the global attractor for 3d viscous primitive equations with viscosity, J Nonlinear Sci, 25(2015), 131--155.



\bibitem {Kob1}
G.M. Kobelkov,  Existence of a solution ¡®in the large¡¯ for the 3D large-scale ocean dynamics
equations, C. R. Math. Acad. Sci. Paris 343(2006), 283--286.

\bibitem {Kob2}
 G.M. Kobelkov,  Existence of a solution ¡®in the large¡¯ for ocean dynamics equations, J. Math. Fluid
Mech. 9(2007), 588--610.

\bibitem {KZ}
I. Kukavica  and M. Ziane,  On the regularity of the primitive equations of the ocean, Nonlinearity
20(2007), 2739--2753.


\bibitem {L}
O. Ladyzhenskaya, ¡°Attractors for Semigroups and Evolution Equations,¡± Cambridge University Press, 1991.

\bibitem {Li}
J. Lions, ¡°Quelques M¡äethode de R¡äesolution des Probl` emes aux Limites Non Lin¡äeaires,¡±
Dunod, Paris, 1969.

\bibitem {LC}
J. Li, J. Chou, Asymptotic behavior of solutions of the moist atmospheric equations, Acta Meteor. Sinica, 56(1998), 61--72
(in Chinese).

\bibitem {LM}
J. Lions and B. Magenes, ¡°Nonhomogeneous Boundary Value Problems and Applications,¡±
Springer--Verlag, New York, 1972.

\bibitem {LTW1}
J.L. Lions, R. Temam, S. Wang, New formulations of the primitive equations of atmosphere and applications, Nonlinearity 5(1992),
 237--288.

\bibitem {LTW2}
J.L. Lions, R. Temam, S. Wang,  On the equations of the large-scale ocean, Nonlinearity
5(1992), 1007--1053.


\bibitem {LTW3}
J.L. Lions, R. Temam, S. Wang, Models of the coupled atmosphere and ocean (CAO I), Comput. Mech. Adv. 1(1993),  1--54.

\bibitem {M}
A. Majda, Introduction to PDEs and Waves for the Atmosphere and Ocean, Courant Lect. Notes Math., vol. 9, 2003.

\bibitem {MT}
A. Miranville  and R. Temam,   Mathematical Modeling in Continuum Mechanics (Cambridge: Cambridge
University Press), 2005.

\bibitem {Pe}
J. Pedlosky, Geophysical Fluid Dynamics, second edition, Springer-Verlag, Berlin/New York, 1987.

\bibitem {PTZ}
 M. Petcu, R. M. Temam  and  M. Ziane,  Some mathematical problems in geophysical fluid dynamics, Handbook
of Numerical Analysis vol 14 Special vol Computational Methods for the Atmosphere and the Oceans
(Amsterdam: Elsevier/North-Holland) pp 577--750, 2009.


\bibitem {Ri}
L.F. Richardson, Weather Prediction by Numerical Process(Cambridge Mathematical Library)(Cambridge:
Cambridge University), 2007.


\bibitem {R}
J. Robinson, ¡°Infinite-dimensional dynamical systems: an introduction to dissipative parabolic PDEs and the theory of global attractors,¡± Cambridge University Press, Cambridge, 2001.


\bibitem {SY}
G. Sell and Y. You, ¡°Dynamics of evolutionary equations,¡± Springer-Verlag, New York, 2002.


\bibitem {T}
R. Temam, ¡° Infinite Dimensional Dynamical Systems in Mechanics and Physics,¡± SpringVerlag, 1988, 2nd Edition, 1997.


\bibitem {TZ}
R. Temam and M. Ziane, Some mathematical problems in geophysical fluid dynamics, Handbook of mathematical fluid dynamics, vol 3, 2004.

\bibitem {W1}
S. Wang, Attractors for the 3-D baroclinic quasi-geostrophic equations of large-scale atmosphere, J. Math. Anal.
Appl. 165 (1) (1992), 266--283.

\bibitem {W2}
J. Wang, Global solutions of the 2D dissipative quasi-geostrophic equations in Besov spaces, SIAM J. Math. Anal. 36 (2004),
1014--1030.

\bibitem {W3}
J. Wang, The two-dimensional quasi-geostrophic equation with critical or supercritical dissipation, Nonlinearity, 18 (2005),
139--154.

\bibitem {ZHKTZ}
M.C. Zelati, A.Huang, L. Kukavica, R.Teman and M. Ziane, The primitive equations of the atmosphere in presence of vapour saturation, Nonlinearlity, 28(2015), 625--668.
\end{thebibliography}
\end{document}